\documentclass[submission,copyright,creativecommons,sharealike]{eptcs}
\providecommand{\event}{ACT 2023} %
\providecommand{\thisvolume}[1]{this volume of EPTCS, Open Publishing Association}
\usepackage{breakurl}             %
\usepackage{underscore}           %

\usepackage[utf8]{inputenc}
\usepackage[T1]{fontenc} %

\usepackage{microtype}
\usepackage{amsfonts}
\usepackage{amssymb}
\usepackage{amsthm}
\usepackage{amsmath}
\usepackage{caption}
\usepackage{etoolbox}
\usepackage{stmaryrd}
\usepackage{ifthen}
\usepackage{suffix}
\usepackage{verbatim} %

\usepackage{xfrac}

\usepackage{bm}
\usepackage{physics}
\usepackage[safe]{tipa} %
\usepackage{mathtools} %
\usepackage{color}
\usepackage{wrapfig}
\usepackage{tikz}
\usepackage{tikzit}
\usetikzlibrary{cd, babel}

\usepackage[framemethod=tikz]{mdframed}

\newcommand*{\secref}[1]{\S\ref{#1}}

\usepackage[backend=biber,style=numeric-comp,sorting=none,date=short]{biblatex}
\appto{\bibsetup}{\sloppy}

\newcommand*{\citep}[1]{\parencite{#1}}
\newcommand*{\citet}[1]{\textcite{#1}}

\usepackage[capitalize]{cleveref} %

\usepackage{graphicx}

\usepackage{pict2e}

\makeatletter
\newcommand{\adjunction}{\@ifstar\named@adjunction\normal@adjunction}
\newcommand{\normal@adjunction}[4]{%
  #1\colon #2%
  \mathrel{\vcenter{%
    \offinterlineskip\m@th
    \ialign{%
      \hfil$##$\hfil\cr
      \longrightharpoonup\cr
      \noalign{\kern-.3ex}
      \smallbot\cr
      \longleftharpoondown\cr
    }%
  }}%
  #3 \noloc #4%
}
\newcommand{\named@adjunction}[4]{%
  #2%
  \mathrel{\vcenter{%
    \offinterlineskip\m@th
    \ialign{%
      \hfil$##$\hfil\cr
      \scriptstyle#1\cr
      \noalign{\kern.1ex}
      \longrightharpoonup\cr
      \noalign{\kern-.3ex}
      \smallbot\cr
      \longleftharpoondown\cr
      \scriptstyle#4\cr
    }%
  }}%
  #3%
}
\newcommand{\longrightharpoonup}{\relbar\joinrel\rightharpoonup}
\newcommand{\longleftharpoondown}{\leftharpoondown\joinrel\relbar}
\newcommand\noloc{%
  \nobreak
  \mspace{6mu plus 1mu}
  {:}
  \nonscript\mkern-\thinmuskip
  \mathpunct{}
  \mspace{2mu}
}
\newcommand{\smallbot}{%
  \begingroup\setlength\unitlength{.15em}%
  \begin{picture}(1,1)
  \roundcap
  \polyline(0,0)(1,0)
  \polyline(0.5,0)(0.5,1)
  \end{picture}%
  \endgroup
}
\makeatother

\let\op=\relax

\def\op{\ensuremath{^{\,\mathrm{op}}}}
\def\coop{\ensuremath{^{\,\mathrm{co\,op}}}}

\newcommand{\rr}{{\mathbb{R}}}

\newcommand{\Cat}[1]{\mathbf{#1}}
\newcommand{\cat}[1]{\mathcal{#1}}
\newcommand{\Fun}[1]{\mathsf{#1}}

\newcommand{\Kl}{\mathcal{K}\mspace{-2mu}\ell}

\DeclareMathOperator*{\E}{\mathbb{E}}

\renewcommand{\d}{\mathrm{d}}

\newcommand{\copara}[1][]{\Cat{Copara}^{#1}}
\newcommand{\ccopara}[1][]{{\Cat{Copara}^{#1}_2}}

\newcommand{\Da}{{\mathcal{D}}}

\DeclareMathOperator{\im}{\mathsf{im}}

\DeclareMathOperator{\id}{\mathsf{id}}

\DeclareMathOperator{\Set}{\Cat{Set}}

\newcommand{\xto}[2][]{\xrightarrow[#1]{#2}}

\makeatletter
\newcommand{\mathoverlap}[2]{\mathpalette\mathoverlap@{{#1}{#2}}}
\newcommand{\mathoverlap@}[2]{\mathoverlap@@{#1}#2}
\newcommand{\mathoverlap@@}[3]{\ooalign{$\m@th#1#2$\crcr\hidewidth$\m@th#1#3$\hidewidth}}
\makeatother

\newcommand{\klcirc}{\bullet} %
\newcommand*{\smallklcirc}{\raisebox{0.18ex}{\scalebox{0.66}{$\klcirc$}}}
\newcommand{\klto}{\mathoverlap{\rightarrow}{\smallklcirc\,}}

\newcommand{\xklto}[2][]{\mathoverlap{\xrightarrow[#1]{#2}}{\smallklcirc\,}}

\newcommand{\lenscirc}{
  \mathbin{\mathoverlap{\circ}{\raisebox{0.375ex}{\scalebox{1.0}[0.33]{$|$}}}}
}
\newcommand{\lensto}{\mathrel{\ooalign{\hfil$\mapstochar\mkern5mu$\hfil\cr$\to$\cr}}}
\newcommand{\xlensto}[2][]{\mathoverlap{\xrightarrow[#1]{#2}}{\raisebox{0.375ex}{\scalebox{1.0}[0.33]{$|$}}\,}}

\makeatletter
\providecommand*{\xmapstofill@}{%
  \arrowfill@{\mapstochar\relbar}\relbar\rightarrow
}
\providecommand*{\xmapsto}[2][]{%
  \ext@arrow 0395\xmapstofill@{#1}{#2}%
}

\makeatletter
\def\slashedarrowfill@#1#2#3#4#5{%
  $\m@th\thickmuskip0mu\medmuskip\thickmuskip\thinmuskip\thickmuskip
   \relax#5#1\mkern-7mu%
   \cleaders\hbox{$#5\mkern-2mu#2\mkern-2mu$}\hfill
   \mathclap{#3}\mathclap{#2}%
   \cleaders\hbox{$#5\mkern-2mu#2\mkern-2mu$}\hfill
   \mkern-7mu#4$%
}
\def\rightslashedarrowfill@{%
  \slashedarrowfill@\relbar\relbar\mapstochar\rightarrow}
\newcommand\xslashedrightarrow[2][]{%
  \ext@arrow 0055{\rightslashedarrowfill@}{#1}{#2}}
\makeatother

\theoremstyle{definition}
\newtheorem{defn}{Definition}[section]
\newtheorem{notation}[defn]{Notation}

\newtheorem{rmk}[defn]{Remark}

\newtheorem*{rmk*}{Remark}

\newtheorem{prop}[defn]{Proposition}
\newtheorem{prop*}{Proposition}

\newtheorem{lemma}[defn]{Lemma}
\newtheorem{thm}[defn]{Theorem}
\newtheorem{cor}[defn]{Corollary}

\newtheorem*{thm*}{Theorem}
\newtheorem*{cor*}{Corollary}

\theoremstyle{remark}

\definecolor{darkblue}{rgb}{0,0,0.7}

\usepackage{bussproofs}
\usepackage{enumitem}
\usepackage{extarrows}

\tikzstyle{xshiftu}=[shift={(#1, 0)}]
\tikzstyle{yshiftu}=[shift={(0, #1)}]
\tikzstyle{dot}=[inner sep=0.0mm, outer sep=0.0mm, minimum size=1mm, draw, shape=circle]
\tikzstyle{copier}=[dot, fill=white, scale=2.0]
\tikzstyle{black copier}=[dot, fill=black, scale=2.0]
\tikzstyle{white dot}=[dot, fill=white, draw=black]
\tikzstyle{action}=[dot, fill=white, scale=0.667, inner sep=0.5mm]
\tikzstyle{box}=[fill=white, draw=black, shape=rectangle]
\tikzstyle{medium box}=[fill=white, draw=black, shape=rectangle, minimum width=1.5cm, minimum height=0.66cm]
\tikzstyle{arrow box}=[fill=white, draw, shape=rectangle, minimum height=5mm, yshift=-0.5mm, minimum width=5mm]
\tikzstyle{effect}=[regular polygon, regular polygon sides=3, draw]
\tikzstyle{state0}=[regular polygon, regular polygon sides=3, draw, shape border rotate=0]
\tikzstyle{state90}=[regular polygon, regular polygon sides=3, draw, shape border rotate=90]
\tikzstyle{state180}=[regular polygon, regular polygon sides=3, draw, shape border rotate=180]
\tikzstyle{state270}=[regular polygon, regular polygon sides=3, draw, shape border rotate=270]
\tikzstyle{scalar}=[diamond, draw, inner sep=1pt]
\tikzstyle{discarder}=[my ground, draw, inner sep=0pt, minimum width=4.2pt, minimum height=11.2pt, anchor=input, rotate=90]
\tikzstyle{discarder0}=[my ground, draw, inner sep=0pt, minimum width=4.2pt, minimum height=11.2pt, anchor=input, rotate=0]

\tikzstyle{pointy1}=[->]
\tikzstyle{midpoint1}=[-, {postaction={decorate,decoration={markings, mark=at position .5 with {\arrow{>}}}}}]
\tikzstyle{midpointy1pointy}=[->, {postaction={decorate,decoration={markings, mark=at position .5 with {\arrow{>}}}}}]
\tikzstyle{dashed1}=[-, dashed]
\tikzstyle{dash-fill1}=[-, dashed, fill={rgb,255: red,157; green,180; blue,180}]
\tikzstyle{dotted1}=[-, dotted]
\tikzstyle{dash-pointy}=[->, dashed]

\input{strings.tikzdefs}

\ifexternalizetikz\tikzexternaldisable\fi
\newsavebox\sbground
\savebox\sbground{%
  \begin{tikzpicture}[baseline=0pt]
    \draw (0,-.1ex) to (0,.85ex)
    node[ground IEC,draw,anchor=input,inner sep=0pt,
    minimum width=3.15pt,minimum height=8.4pt,rotate=90] {};
  \end{tikzpicture}%
}
\newcommand{\ground}{\mathord{\usebox\sbground}}
\newcommand{\smallground}{\scalebox{0.5}{$\mathord{\usebox\sbground}$}}

\newsavebox\sbcopier
\savebox\sbcopier{%
  \begin{tikzpicture}[baseline=0pt]
    \node[copier,scale=0.7] (a) at (0,3.8pt) {};
    \draw (a) -- +(-90:.21);
    \draw (a) -- +(45:.21);
    \draw (a) -- +(135:.21);
  \end{tikzpicture}}

\ifexternalizetikz\tikzexternalenable\fi

\newsavebox\bsbcopier
\savebox\bsbcopier{%
  \begin{tikzpicture}[baseline=0pt]
    \node[black copier,scale=0.7] (a) at (0,3.8pt) {};
    \draw (a) -- +(-90:.21);
    \draw (a) -- +(45:.21);
    \draw (a) -- +(135:.21);
  \end{tikzpicture}}
\newcommand{\bcopier}{\mathord{\usebox\bsbcopier}}
\ifexternalizetikz\tikzexternalenable\fi

\usepackage{mathabx}
\usepackage{quiver}

\def\proj{\mathsf{proj}}

\def\Stat{\mathsf{Stat}}
\def\SStat{{\Stat_2}}
\def\disc{\Fun{disc}\,}

\def\BLens{\Cat{BayesLens}}

\def\Loss{\mathsf{Loss}}
\def\MonLoss{\mathsf{MonLoss}}
\def\pLoss{\pi_{\Loss}}
\def\pLens{\pi_{\mathsf{Lens}}}

\newcommand{\SGame}[1]{\Cat{SGame}_{#1}}

\def\KL{\mathsf{KL}}
\def\MLE{\mathsf{MLE}}
\def\FE{\mathsf{FE}}
\def\LFE{\mathsf{LFE}}

\author{Toby St Clere Smithe
  \institute{VERSES Research}
  \institute{Topos Institute}
  \email{q@tsmithe.net}
}
\def\authorrunning{T. St Clere Smithe}

\hypersetup{
  pdftitle={Statistical Games},
  pdfauthor={Toby St Clere Smithe}
}

\date{29 June 2023}

\title{Approximate Inference via Fibrations of Statistical Games}
\def\titlerunning{Statistical Games}

\begin{document}

\maketitle

\begin{abstract}
  We characterize a number of well known systems of approximate inference as \textit{loss models}: lax sections of 2-fibrations of statistical games, constructed by attaching internally-defined loss functions to Bayesian lenses.
  Our examples include the relative entropy, which constitutes a \textit{strict} section, and whose chain rule is formalized by the horizontal composition of the 2-fibration.
  In order to capture this compositional structure, we first introduce the notion of `copy-composition', alongside corresponding bicategories through which the composition of copy-discard categories factorizes.
  These bicategories are a variant of the $\copara$ construction, and so we additionally introduce coparameterized Bayesian lenses, proving that coparameterized Bayesian updates compose optically, as in the non-coparameterized case.
\end{abstract}

\section{Introduction}

In previous work \parencite{Smithe2020Cyber}, we introduced \textit{Bayesian lenses}, observing that the Bayesian inversion of a composite stochastic channel is (almost surely) equal to the `lens composite' of the inversions of the factors; that is, \textit{Bayesian updates compose optically} (`BUCO') \parencite{Smithe2020Bayesian}.
Formalizing this statement for a given category $\cat{C}$ all of whose morphisms (`channels') admit Bayesian inversion, we can observe that there is (almost surely) a functor $(-)^\dag:\cat{C}\to\BLens(\cat{C})$ from $\cat{C}$ to the category $\BLens(\cat{C})$ whose morphisms $(X,A)\lensto(Y,B)$ are Bayesian lenses:
pairs $(c,c')$ of a channel $X\klto Y$ with a `state-dependent' inverse $c':\cat{C}(I,X)\to\cat{C}(B,A)$.
Bayesian lenses constitute the morphisms of a fibration $\pLens:\BLens(\cat{C})\to\cat{C}$, since $\BLens(\cat{C})$ is obtained as the Grothendieck construction of (the pointwise opposite of) an indexed category $\Stat:\cat{C}\op\to\Cat{Cat}$ of `state-dependent channels' (recalled in Appendix \ref{sec:stat}), and the functor $(-)^\dag$ is in fact a section of $\pLens$, taking $c:X\klto Y$ to the lens $(c,c^\dag):(X,X)\lensto(Y,Y)$, where $c^\dag$ is the almost-surely unique Bayesian inversion of $c$ (so that the projection $\pLens$ can simply forget the inversion, returning again the channel $c$).

The functor $(-)^\dag$ picks out a special class of Bayesian lenses, which we may call \textit{exact} (as they compute `exact' inversions), but although the category $\BLens(\cat{C})$ has many other morphisms, the construction is not extravagant:
by comparison, we can think of the non-exact lenses as representing \textit{approximate} inference systems.
This is particularly necessary in computational applications, because computing exact inversions is usually intractable, but this creates a new problem: choosing an approximation, and measuring its performance.
In this paper, we formalize this process, by attaching \textit{loss functions} to Bayesian lenses, thus creating another fibration, of \textit{statistical games}.
Sections of this latter fibration encode compositionally well-behaved systems of approximation that we call \textit{loss models}.

A classic example of a loss model will be supplied by the relative entropy, which in some sense measures the `divergence' between distributions: the game here is then to minimize the divergence between the approximate and exact inversions.
If $\pi$ and $\pi'$ are two distributions on a space $X$, with corresponding density functions $p_\pi$ and $p_{\pi'}$ (both with respect to a common measure),
then their relative entropy $D(\pi,\pi')$ is the real number given by $\E_{x\sim\pi}\left[\log p_\pi(x) - \log p_{\pi'}(x)\right]$\footnote{
  For details about this `expectation' notation $\E$, see \ref{nota:expectations}.}.
Given a pair of channels $\alpha,\alpha':A\klto B$ (again commensurately associated with densities), we can extend $D$ to a map $D_{\alpha,\alpha'}:A\to\rr_+$ in the natural way, writing $a\mapsto D\bigl(\alpha(a),\alpha'(a)\bigr)$.
We can assign such a map $D_{\alpha,\alpha'}$ to any such parallel pair of channels, and so, following the logic of composition in $\cat{C}$, we might hope for the following equation to hold for all $a:A$ and composable parallel pairs $\alpha,\alpha':A\klto B$ and $\beta,\beta':B\klto C$,:
\[ D_{\beta\klcirc\alpha,\beta'\klcirc\alpha'}(a) = \E_{b\sim\alpha(a)} \left[ D_{\beta,\beta'}(b) \right] + D_{\alpha,\alpha'}(a) \]

The right-hand side is known as the \textit{chain rule} for relative entropy, but, unfortunately, the equation does \textit{not} hold in general, because the composites $\beta\klcirc\alpha$ and $\beta'\klcirc\alpha'$ involve an extra expectation (by the `Chapman-Kolmogorov' rule for channel composition).
However, we \textit{can} satisfy an equation of this form by using `copy-composition':
if we write $\bcopier_B$ to denote the canonical `copying' comultiplication on $B$, and define $\beta\klcirc^2\alpha := (\id_B\otimes\beta)\klcirc\bcopier_B\klcirc\alpha$, then $D_{\beta\klcirc^2\alpha,\beta'\klcirc^2\alpha'}(a)$ \textit{does} equal the chain-rule form on the right-hand side.
This result exhibits a general pattern about ``copy-discard categories'' \parencite{Cho2017Disintegration} such as $\cat{C}$: composition $\klcirc$ can be decomposed into first copying $\bcopier$, and then discarding $\ground$.
If we don't discard, then we retain the `intermediate' variables, and this results in a functorial assignment of relative entropies to channels.

The rest of this paper is dedicated to making use of this observation to construct loss models, including (but not restricted to) the relative entropy.
The first complication that we encounter is that copy-composition is not strictly unital, because composing with an identity retains an extra variable.
To deal with this, in \secref{sec:copy-comp}, we construct a \textit{bicategory} of copy-composite channels, extending the $\copara$ construction \parencite[\S2]{Capucci2022Foundations}, and build coparameterized (copy-composite) Bayesian lenses accordingly; we also prove a corresponding BUCO result.
In \secref{sec:stat-games}, we then construct 2-fibrations of statistical games, defining loss functions internally to any copy-discard category $\cat{C}$ that admits ``bilinear effects''.
Because we are dealing with approximate systems, the 2-dimensional structure of the construction is useful: loss models are allowed to be \textit{lax} sections.
We then characterize the relative entropy, maximum likelihood estimation, the free energy, and the `Laplacian' free energy as such loss models.

Assuming $\cat{C}$ is symmetric monoidal, the constructions here result in monoidal (2-)fibrations, but due to space constraints we defer the presentation of this structure (and its exemplification by the foregoing loss models) to Appendix \ref{sec:mon-stat-games}.

\begin{rmk}
  Much of this work is situated amongst monoidal fibrations of bicategories, the full theory of which is not known to the present author.
  Fortunately, enough structure is known for the present work to have been possible, and where things become murkier---such as in the context of monoidal indexed bicategories and their lax homomorphisms---the way largely seems clear.
  For this, we are grateful to \textcite{Bakovic2010Fibrations}, \textcite{Johnson20202Dimensional}, and \textcite{Moeller2018Monoidal} in particular for lighting the way; and we enthusiastically encourage the further elucidation of these structures by category theorists.
\end{rmk}

\begin{rmk}
  For reasons of space, detailed proofs are not included in the proceedings version of this paper; however, they are included in an appendix to the conference submission, which is available on the arXiv repository with the paper ID \href{https://arxiv.org/abs/2306.17009v1}{2306.17009v1}.
\end{rmk}

\section{`Copy-composite' Bayesian lenses} \label{sec:copy-comp}

\subsection{Copy-composition by coparameterization}

In a locally small copy-discard category $\cat{C}$, every object $A$ is equipped with a canonical comonoid structure $(\bcopier_A,\ground_A)$, and so, by the comonoid laws, it is almost a triviality that the composition function $\klcirc:\cat{C}(B,C)\times\cat{C}(A,B)\to\cat{C}(A,C)$ factorizes as
\begin{gather*}
  \cat{C}(B,C)\times\cat{C}(A,B)
  \xto{(\id_B\otimes-)\times\cat{C}\left(\id_A,\bcopier_B\right)}
  \cat{C}(B\otimes B,B\otimes C)\times\cat{C}(A,B\otimes B)
  \; \cdots \\ \cdots \; \xto{\klcirc}
  \cat{C}(A,B\otimes C)
  \xto{\cat{C}(\id_A,\proj_C)}
  \cat{C}(A,C)
\end{gather*}
where the first factor copies the $B$ output of the first morphism and tensors the second morphism with the identity on $B$, the second factor composes the latter tensor with the copies, and the third discards the extra copy of $B$\footnote{\label{fn:proj} We define $\proj_C := B\otimes C \xto{\ground_B\otimes\id_C} I\otimes C \xto{\lambda_C} C$, using the comonoid counit and the left unitor of $\cat{C}$'s monoidal structure.}.
This is, however, only \textit{almost} trivial, since it witnesses the structure of `Chapman-Kolmogorov' style composition in categories of stochastic channels such as $\Kl(\Da)$, the Kleisli category of the (finitary) distributions monad $\Da:\Cat{Set}\to\Cat{Set}$.
There, given channels $c:A\klto B$ and $d:B\klto C$, the composite $d\klcirc c$ is formed first by constructing the `joint' channel $d\klcirc^2c$ defined by $(d\klcirc^2c)(b,c|a) := d(c|b)c(b|a)$, and then discarding (marginalizing over) $b:B$, giving
\[ (d\klcirc c)(c|a) = \sum_{b:B} (d\klcirc^2c)(b,c|a) = \sum_{b:B} d(c|b)c(b|a) \, . \]
Of course, the channel $d\klcirc^2 c$ is not a morphism $A\klto C$, but rather $A\klto B\otimes C$; that is, it is \textit{coparameterized} by $B$.
Moreover, as noted above, $\klcirc^2$ is not strictly unital: we need a 2-cell that discards the coparameter, and hence a bicategory, in order to recover (weak) unitality.
We therefore construct a bicategory $\ccopara(\cat{C})$ as a variant of the $\copara$ construction \parencite[\S2]{Capucci2022Foundations}, in which a 1-cell $A\to B$ may be any morphism $A\klto M\otimes B$ in $\cat{C}$, and where horizontal composition is precisely copy-composition.

\begin{thm} \label{thm:copara2}
  Let $(\cat{C},\otimes,I)$ be a copy-discard category.
  Then there is a bicategory $\ccopara(\cat{C})$ as follows.
  Its 0-cells are the objects of $\cat{C}$.
  A 1-cell $f:A\xto[M]{}B$ is a morphism $f:A\to M\otimes B$ in $\cat{C}$.
  A 2-cell $\varphi:f\Rightarrow f'$, with $f:A\xto[M]{}B$ and $f':A\xto[M']{}B$, is a morphism $\varphi:A\otimes M\otimes B\to M'$ of $\cat{C}$, satisfying the \textit{change of coparameter} axiom:
  \[ \scalebox{0.75}
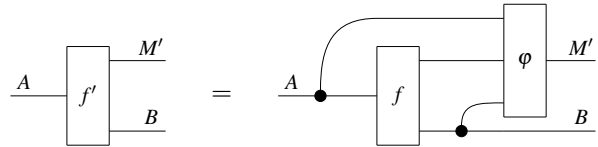
{\begin{tikzpicture}
	\begin{pgfonlayer}{nodelayer}
		\node [style=none] (0) at (0.75, 1.25) {};
		\node [style=none] (1) at (0.75, 1.75) {};
		\node [style=none] (2) at (-0.75, 1.75) {};
		\node [style=none] (3) at (-0.75, 0) {};
		\node [style=none] (4) at (-0.75, -1.75) {};
		\node [style=none] (5) at (0.75, -1.75) {};
		\node [style=none] (6) at (0.75, -1.25) {};
		\node [style=none] (7) at (0, 0) {$f'$};
		\node [style=none] (8) at (-2.75, 0) {};
		\node [style=none] (13) at (-2.25, 0.5) {$A$};
		\node [style=none] (24) at (2.75, -1.25) {};
		\node [style=none] (25) at (2.25, -0.75) {$B$};
		\node [style=none] (27) at (2.75, 1.25) {};
		\node [style=none] (28) at (2.25, 1.75) {$M'$};
	\end{pgfonlayer}
	\begin{pgfonlayer}{edgelayer}
		\draw (8.center) to (3.center);
		\draw (2.center) to (1.center);
		\draw (4.center) to (5.center);
		\draw (2.center) to (4.center);
		\draw (1.center) to (5.center);
		\draw (0.center) to (27.center);
		\draw (6.center) to (24.center);
	\end{pgfonlayer}
\end{tikzpicture}} \quad\,=\quad\, \scalebox{0.75}{\begin{tikzpicture}
	\begin{pgfonlayer}{nodelayer}
		\node [style=none] (0) at (0.75, 1.25) {};
		\node [style=none] (1) at (0.75, 1.75) {};
		\node [style=none] (2) at (-0.75, 1.75) {};
		\node [style=none] (3) at (-0.75, 0) {};
		\node [style=none] (4) at (-0.75, -1.75) {};
		\node [style=none] (5) at (0.75, -1.75) {};
		\node [style=none] (6) at (0.75, -1.25) {};
		\node [style=none] (7) at (0, 0) {$f$};
		\node [style=none] (13) at (-3.75, 0.5) {$A$};
		\node [style=black copier] (24) at (2.25, -1.25) {};
		\node [style=none] (25) at (6.5, -0.75) {$B$};
		\node [style=none] (27) at (3.75, 1.25) {};
		\node [style=none] (28) at (6.5, 1.75) {$M'$};
		\node [style=none] (29) at (5.25, 1.25) {};
		\node [style=none] (35) at (3.75, -1.25) {};
		\node [style=none] (36) at (4.5, 1.25) {$\varphi$};
		\node [style=none] (37) at (7, 1.25) {};
		\node [style=none] (38) at (7, -1.25) {};
		\node [style=none] (39) at (3.75, 2.75) {};
		\node [style=black copier] (46) at (-2.75, 0) {};
		\node [style=none] (47) at (-4.25, 0) {};
		\node [style=none] (54) at (3.75, -0.25) {};
		\node [style=none] (55) at (3.75, 2.75) {};
		\node [style=none] (56) at (3.75, 3.25) {};
		\node [style=none] (57) at (3.75, -0.75) {};
		\node [style=none] (58) at (5.25, 3.25) {};
		\node [style=none] (59) at (5.25, -0.75) {};
		\node [style=none] (60) at (-0.75, 2.75) {};
	\end{pgfonlayer}
	\begin{pgfonlayer}{edgelayer}
		\draw (2.center) to (1.center);
		\draw (4.center) to (5.center);
		\draw (2.center) to (4.center);
		\draw (1.center) to (5.center);
		\draw (0.center) to (27.center);
		\draw (6.center) to (24);
		\draw (29.center) to (37.center);
		\draw (35.center) to (38.center);
		\draw (47.center) to (46);
		\draw (46) to (3.center);
		\draw (57.center) to (56.center);
		\draw [in=180, out=90, looseness=1.25] (24) to (54.center);
		\draw (24) to (35.center);
		\draw (58.center) to (59.center);
		\draw (57.center) to (59.center);
		\draw (56.center) to (58.center);
		\draw (46)
			 to [in=-180, out=90, looseness=1.25] (60.center)
			 to (55.center);
	\end{pgfonlayer}
\end{tikzpicture}} \]
  The identity 2-cell $\id_f:f\Rightarrow f$ on $f:A\xto[M]{}B$ is given by the projection morphism $\proj_M:A\otimes M\otimes B\to M$ obtained by discarding $A$ and $B$, as in footnote \ref{fn:proj}.
  The identity 1-cell $\id_A$ on $A$ is given by the inverse of the left unitor of the monoidal structure on $\cat{C}$, \textit{i.e.} $\id_A := \lambda_A^{-1} : A\xto[I]{}A$, with coparameter thus given by the unit object $I$.

  Given 2-cells $\varphi:f\Rightarrow f'$ and $\varphi':f'\Rightarrow f''$, their vertical composite $\varphi'\odot\varphi:f\Rightarrow f''$ is given by the string diagram on the left below.
  Given 1-cells $f:A\xto[M]{}B$ then $g:B\xto[N]{}C$, the horizontal composite $g\circ f:A\xto[(M\otimes B)\otimes N]{}C$ is given by the middle string diagram below.
  Given 2-cells $\varphi:f\Rightarrow f'$ and $\gamma:g\Rightarrow g'$ between 1-cells $f,f':A\xto[M]{}B$ and $g,g':B\xto[N]{}C$, their horizontal composite $\gamma\circ\varphi:(g\circ f)\Rightarrow (g'\circ f')$ is defined by the string diagram on the right below.
  \[ \scalebox{0.75}{\begin{tikzpicture}
	\begin{pgfonlayer}{nodelayer}
		\node [style=none] (0) at (-0.75, 0) {};
		\node [style=none] (1) at (0.75, 0) {};
		\node [style=none] (2) at (0, 0) {$\varphi$};
		\node [style=none] (4) at (-0.75, -1.25) {};
		\node [style=none] (6) at (-0.75, 1.75) {};
		\node [style=none] (7) at (-0.75, -1.75) {};
		\node [style=none] (8) at (0.75, 1.75) {};
		\node [style=none] (9) at (0.75, -1.75) {};
		\node [style=none] (10) at (2.25, 0) {};
		\node [style=none] (11) at (3.75, 0) {};
		\node [style=none] (12) at (3, 0) {$\varphi'$};
		\node [style=none] (13) at (2.25, -2.75) {};
		\node [style=none] (14) at (2.25, 2.75) {};
		\node [style=none] (15) at (2.25, 3.25) {};
		\node [style=none] (16) at (2.25, -3.25) {};
		\node [style=none] (17) at (3.75, 3.25) {};
		\node [style=none] (18) at (3.75, -3.25) {};
		\node [style=black copier] (19) at (-2.25, -2) {};
		\node [style=none] (20) at (-0.75, -2.75) {};
		\node [style=none] (21) at (-3.75, -2) {};
		\node [style=none] (22) at (-3.75, 0) {};
		\node [style=none] (25) at (5.25, 0) {};
		\node [style=none] (26) at (-0.75, 2.75) {};
		\node [style=black copier] (27) at (-2.25, 2) {};
		\node [style=none] (28) at (-0.75, 1.25) {};
		\node [style=none] (29) at (-3.75, 2) {};
	\end{pgfonlayer}
	\begin{pgfonlayer}{edgelayer}
		\draw (7.center) to (6.center);
		\draw (8.center) to (9.center);
		\draw (7.center) to (9.center);
		\draw (6.center) to (8.center);
		\draw (16.center) to (15.center);
		\draw (17.center) to (18.center);
		\draw (16.center) to (18.center);
		\draw (15.center) to (17.center);
		\draw (4.center)
			 to [in=90, out=-180, looseness=1.25] (19.center)
			 to [in=180, out=-90, looseness=1.25] (20.center)
			 to (13.center);
		\draw (21.center) to (19);
		\draw (22.center) to (0.center);
		\draw (1.center) to (10.center);
		\draw (11.center) to (25.center);
		\draw (28.center)
			 to [in=-90, out=-180, looseness=1.25] (27.center)
			 to [in=-180, out=90, looseness=1.25] (26.center)
			 to (14.center);
		\draw (29.center) to (27);
	\end{pgfonlayer}
\end{tikzpicture}} \hspace*{1.5cm} \scalebox{0.75}{\begin{tikzpicture}
	\begin{pgfonlayer}{nodelayer}
		\node [style=none] (0) at (0.75, 1) {};
		\node [style=none] (1) at (0.75, 1.5) {};
		\node [style=none] (2) at (-0.75, 1.5) {};
		\node [style=none] (3) at (-0.75, -0.25) {};
		\node [style=none] (4) at (-0.75, -2) {};
		\node [style=none] (5) at (0.75, -2) {};
		\node [style=none] (6) at (0.75, -1.5) {};
		\node [style=none] (7) at (0, -0.25) {$f$};
		\node [style=none] (8) at (-2.75, -0.25) {};
		\node [style=none] (9) at (3.25, 3.5) {};
		\node [style=none] (13) at (-2.25, 0.25) {$A$};
		\node [style=none] (14) at (4.75, 1) {};
		\node [style=none] (15) at (4.75, 1.5) {};
		\node [style=none] (16) at (3.25, 1.5) {};
		\node [style=none] (17) at (3.25, -0.25) {};
		\node [style=none] (18) at (3.25, -2) {};
		\node [style=none] (19) at (4.75, -2) {};
		\node [style=none] (20) at (4.75, -1.5) {};
		\node [style=none] (21) at (4, -0.25) {$g$};
		\node [style=none] (23) at (6.75, 1) {};
		\node [style=none] (24) at (6.75, -1.5) {};
		\node [style=none] (25) at (6.25, -1) {$C$};
		\node [style=none] (26) at (6.25, 1.5) {$N$};
		\node [style=none] (27) at (6.75, 3.5) {};
		\node [style=none] (28) at (6.25, 4) {$M$};
		\node [style=black copier] (29) at (2, -1) {};
		\node [style=none] (30) at (6.75, 2.25) {};
		\node [style=none] (31) at (6.25, 2.75) {$B$};
		\node [style=none] (32) at (4, 2.25) {};
	\end{pgfonlayer}
	\begin{pgfonlayer}{edgelayer}
		\draw (8.center) to (3.center);
		\draw (2.center) to (1.center);
		\draw (4.center) to (5.center);
		\draw (2.center) to (4.center);
		\draw (1.center) to (5.center);
		\draw [in=180, out=0, looseness=1.25] (0.center) to (9.center);
		\draw (20.center) to (24.center);
		\draw (16.center) to (15.center);
		\draw (18.center) to (19.center);
		\draw (16.center) to (18.center);
		\draw (15.center) to (19.center);
		\draw (14.center) to (23.center);
		\draw (9.center) to (27.center);
		\draw [in=-180, out=90, looseness=1.25] (29) to (32.center);
		\draw (32.center) to (30.center);
		\draw [in=-150, out=0, looseness=1.25] (6.center) to (29);
		\draw [in=-180, out=45, looseness=1.25] (29) to (17.center);
	\end{pgfonlayer}
\end{tikzpicture}} \hspace*{1.5cm} \scalebox{0.75}{\begin{tikzpicture}
	\begin{pgfonlayer}{nodelayer}
		\node [style=none] (0) at (2, 2.5) {};
		\node [style=none] (1) at (3.5, 2.5) {};
		\node [style=none] (2) at (2.75, 2.5) {$\varphi$};
		\node [style=none] (3) at (0.5, 1.25) {};
		\node [style=none] (4) at (2, 3.75) {};
		\node [style=none] (5) at (2, 4.25) {};
		\node [style=none] (6) at (2, 0.75) {};
		\node [style=none] (7) at (3.5, 4.25) {};
		\node [style=none] (8) at (3.5, 0.75) {};
		\node [style=none] (16) at (2, -2.5) {};
		\node [style=none] (17) at (3.5, -2.5) {};
		\node [style=none] (18) at (2.75, -2.5) {$\gamma$};
		\node [style=none] (19) at (2, -3.75) {};
		\node [style=none] (20) at (2, -1.25) {};
		\node [style=none] (21) at (2, -0.75) {};
		\node [style=none] (22) at (2, -4.25) {};
		\node [style=none] (23) at (3.5, -0.75) {};
		\node [style=none] (24) at (3.5, -4.25) {};
		\node [style=none] (25) at (2, 0) {};
		\node [style=black copier] (26) at (0.5, 0) {};
		\node [style=black copier] (27) at (-1, 0) {};
		\node [style=none] (28) at (-2.5, 0) {};
		\node [style=none] (29) at (-2.5, -2.5) {};
		\node [style=none] (30) at (-2.5, -3.75) {};
		\node [style=none] (31) at (-2.5, 2.5) {};
		\node [style=none] (32) at (-2.5, 3.75) {};
		\node [style=none] (33) at (2, 1.25) {};
		\node [style=none] (34) at (5, 2.5) {};
		\node [style=none] (35) at (5, 0) {};
		\node [style=none] (36) at (5, -2.5) {};
		\node [style=none] (37) at (-2, 4.25) {$A$};
		\node [style=none] (38) at (-2, 3) {$M$};
		\node [style=none] (39) at (-2, 0.5) {$B$};
		\node [style=none] (40) at (-2, -2) {$N$};
		\node [style=none] (41) at (-2, -3.25) {$C$};
		\node [style=none] (42) at (4.5, -2) {$N'$};
		\node [style=none] (43) at (4.5, 0.5) {$B$};
		\node [style=none] (44) at (4.5, 3) {$M'$};
	\end{pgfonlayer}
	\begin{pgfonlayer}{edgelayer}
		\draw (6.center) to (5.center);
		\draw (7.center) to (8.center);
		\draw (6.center) to (8.center);
		\draw (5.center) to (7.center);
		\draw (22.center) to (21.center);
		\draw (23.center) to (24.center);
		\draw (22.center) to (24.center);
		\draw (21.center) to (23.center);
		\draw (28.center) to (27);
		\draw (27)
			 to [in=180, out=90, looseness=1.25] (3.center)
			 to (33.center);
		\draw (27) to (26);
		\draw (26) to (25.center);
		\draw (32.center) to (4.center);
		\draw (31.center) to (0.center);
		\draw (29.center) to (16.center);
		\draw (30.center) to (19.center);
		\draw (17.center) to (36.center);
		\draw (25.center) to (35.center);
		\draw (1.center) to (34.center);
		\draw [in=180, out=-90, looseness=1.25] (26) to (20.center);
	\end{pgfonlayer}
\end{tikzpicture}} \]
\end{thm}

\begin{rmk}
  When $\cat{C}$ is symmetric monoidal, $\ccopara(\cat{C})$ inherits a monoidal structure, elaborated in Proposition \ref{prop:mon-ccopara}.
\end{rmk}

\begin{rmk}
  In order to capture the bidirectionality of Bayesian inversion we will need to conside left- and right-handed versions of the $\ccopara$ construction.
  These are formally dual, and when $\cat{C}$ is symmetric monoidal (as in most examples) they are isomorphic.
  Nonetheless, it makes formalization easier if we explicitly distinguish $\ccopara[l](\cat{C})$, in which the coparameter is placed on the left of the codomain (as above), from $\ccopara[r](\cat{C})$, in which it is placed on the right.
  Aside from the swapping of this handedness, the rest of the construction is the same.
\end{rmk}

We end this section with three easy (and ambidextrous) propositions relating $\cat{C}$ and $\ccopara(\cat{C})$.

\begin{prop}
  There is an identity-on-objects lax embedding $\iota:\cat{C}\hookrightarrow\ccopara(\cat{C})$, mapping $f:X\to Y$ to $f:X\xto[I]{}Y$ (using the unitor of the monoidal structure on $\cat{C}$).
  The laxator $\iota(g)\circ\iota(f)\to\iota(g\circ f)$ discards the coparameter obtained from copy-composition.
\end{prop}

\begin{prop} \label{prop:discard-func}
  There is a `discarding' functor $(-)^{\smallground}:\ccopara(\cat{C})\to\cat{C}$, which takes any coparameterized morphism and discards the coparameter.
\end{prop}

\begin{prop}
  $\iota$ is a section of $(-)^{\smallground}$. That is, $\id_{\cat{C}} = \cat{C}\xhookrightarrow{\iota}\ccopara(\cat{C})\xto{(-)^{\smallground}}\cat{C}$.
\end{prop}

\subsection{Coparameterized Bayesian lenses}

In order to define (bi)categories of statistical games, coherently with loss functions like the relative entropy that compose by copy-composition, we first need to define coparameterized (copy-composite) Bayesian lenses.
Analogously to non-coparameterized Bayesian lenses, these will be obtained by applying a Grothendieck construction to an indexed bicategory \parencite[{Def. 3.5}]{Bakovic2010Fibrations} of state-dependent channels.

\begin{defn}
  We define the indexed bicategory $\SStat:\ccopara[l](\cat{C})\coop\to\Cat{Bicat}$ fibrewise as follows.
  \begin{enumerate}[label=(\roman*)]
  \item The 0-cells $\SStat(X)_0 $ of each fibre $\SStat(X)$ are the objects $\cat{C}_0$ of $\cat{C}$.
  \item For each pair of 0-cells $A,B$, the hom-category $\SStat(X)(A,B)$ is defined to be the functor category $\Cat{Cat}\bigl(\disc\cat{C}(I,X),\ccopara[r](\cat{C})(A,B)\bigr)$, where $\disc$ denotes the functor taking a set to the associated discrete category.
  \item For each 0-cell $A$, the identity functor $\id_A : \Cat{1}\to\SStat(X)(A,A)$ is the constant functor on the identity on $A$ in $\ccopara[r](\cat{C})$; \textit{i.e.} $\disc\cat{C}(I,X) \xto{!} 1 \xto{\id_A} \ccopara[r](\cat{C})(A,A)$.
  \item For each triple $A,B,C$ of 0-cells, the horizontal composition functor $\circ_{A,B,C}$ is defined by
    \begin{gather*}
      \circ_{A,B,C} : \SStat(X)(B,C)\times\SStat(X)(A,B) \;\cdots \\
      \cdots\; \xto{=} \Cat{Cat}\bigl(\disc\cat{C}(I,X),\ccopara[r](\cat{C})(B,C)\bigr) \times \Cat{Cat}\bigl(\disc\cat{C}(I,X),\ccopara[r](\cat{C})(A,B)\bigr) \;\cdots \\
      \cdots\; \xto{\times} \Cat{Cat}\bigl(\disc\cat{C}(I,X)^2,\ccopara[r](\cat{C})(B,C)\times\ccopara[r](\cat{C})(A,B)\bigr) \;\cdots \\
      \cdots\; \xto{\Cat{Cat}\left(\bcopier,\circ\right)} \Cat{Cat}\bigl(\disc\cat{C}(I,X),\ccopara[r](\cat{C})(A,C)\bigr) \;\cdots \\
      \cdots\; \xto{=} \SStat(X)(A,C)
    \end{gather*}
    where $\Cat{Cat}\left(\bcopier,\circ\right)$ indicates pre-composition with the universal (Cartesian) copying functor in $(\Cat{Cat},\times,\Cat{1})$ and post-composition with the horizontal composition functor of $\ccopara[r](\cat{C})$.
  \end{enumerate}
  For each pair of 0-cells $X,Y$ in $\copara[l](\cat{C})$, we define the reindexing pseudofunctor $\SStat_{,X,Y}:\copara[l](\cat{C})(X,Y)\op\to\Cat{Bicat}\bigl(\SStat(Y),\SStat(X)\bigr)$ as follows.
  \begin{enumerate}[label=(\alph*)]
  \item For each 1-cell $f$ in $\copara[l](\cat{C})(X,Y)$, we obtain a pseudofunctor $\SStat(f):\SStat(Y)\to\SStat(X)$ which acts as the identity on 0-cells.
  \item For each pair of 0-cells $A,B$ in $\SStat(Y)$, the functor $\SStat(f)_{A,B}$ is defined as the precomposition functor $\Cat{Cat}\bigl(\disc\cat{C}(I,f^{\smallground}),\ccopara[r](\cat{C})(A,B)\bigr)$, where $(-)^{\smallground}$ is the discarding functor $\ccopara[l](\cat{C})\to\cat{C}$ of Proposition \ref{prop:discard-func}.
    \item For each 2-cell $\varphi:f\Rightarrow f'$ in $\ccopara[l](\cat{C})(X,Y)$, the pseudonatural transformation $\SStat(\varphi):\SStat(f')\Rightarrow\SStat(f)$ is defined on 0-cells $A:\SStat(Y)$ by the discrete natural transformation with components $\SStat(\varphi)_A := \id_A$, and on 1-cells $c:\SStat(Y)(A,B)$ by the substitution natural transformation with constitutent 2-cells $\SStat(\varphi)_c:\SStat(f)(c)\Rightarrow\SStat(f')(c)$ in $\SStat(X)$ which acts by replacing $\Cat{Cat}\bigl(\disc\cat{C}(I,f^{\smallground}),\ccopara[r](\cat{C})(A,B)\bigr)$ by $\Cat{Cat}\bigl(\disc\cat{C}(I,f'^{\smallground}),\ccopara[r](\cat{C})(A,B)\bigr)$; and which we might alternatively denote by $\Cat{Cat}\bigl(\disc\cat{C}(I,\varphi^{\smallground}),\ccopara[r](\cat{C})(A,B)\bigr)$.
  \end{enumerate}
\end{defn}

\begin{notation}
  We will write $f:A\xklto[M]{X}B$ to denote a state-dependent coparameterized channel $f$ with coparameter $M$ and state-dependence on $X$.
\end{notation}

In 1-category theory, lenses are morphisms in the fibrewise opposite of a fibration \parencite{Spivak2019Generalized}.
Analogously, our bicategorical Bayesian lenses are obtained as 1-cells in the bicategorical Grothendieck construction \parencite[\S6]{Bakovic2010Fibrations} of (the pointwise opposite of) the indexed bicategory $\SStat$.

\begin{defn}
  Fix a copy-discard category $(\cat{C},\otimes,I)$.
  We define the bicategory of coparameterized Bayesian lenses in $\cat{C}$, denoted $\BLens_2(\cat{C})$ or simply $\BLens_2$, to be the bicategorical Grothendieck construction of the pointwise opposite of the corresponding indexed bicategory $\SStat$, with the following data:
  \begin{enumerate}[label=(\roman*)]
  \item A 0-cell in $\BLens_2$ is a pair $(X,A)$ of an object $X$ in $\ccopara[l](\cat{C})$ and an object $A$ in $\SStat(X)$; equivalently, a 0-cell in $\BLens_2$ is a pair of objects in $\cat{C}$.
  \item The hom-category $\BLens_2\bigl((X,A),(Y,B)\bigr)$ is the product category $\ccopara[l](\cat{C})(X,Y)\times\SStat(X)(B,A)$.
  \item The identity on $(X,A)$ is given by the pair $(\id_X,\id_A)$.
  \item For each triple of 0-cells $(X,A),(Y,B),(Z,C)$, the horizontal composition functor is given by
    \begin{align*}
      & \BLens_2\bigl((Y,B),(Z,C)\bigr)\times\BLens_2\bigl((X,A),(Y,B)\bigr) \\
      =\quad& \ccopara[l](\cat{C})(Y,Z)\times\SStat(Y)(C,B) \times \ccopara[l](\cat{C})(X,Y)\times\SStat(X)(B,A) \\
      \xto{\sim}\quad& \sum_{g:\ccopara[l](\cat{C})(Y,Z)} \sum_{f:\ccopara[l](\cat{C})(X,Y)} \SStat(Y)(C,B) \times \SStat(X)(B,A) \\
      \xto{\sum_g \sum_f \SStat(f)_{C,B}\times\id}\quad& \sum_{g:\ccopara[l](\cat{C})(Y,Z)} \sum_{f:\ccopara[l](\cat{C})(X,Y)} \SStat(X)(C,B) \times \SStat(X)(B,A) \\
      \xto{\sum_{\circ^{\ccopara[l](\cat{C})}} \circ^{\SStat(X)}}\quad& \sum_{g\circ f:\ccopara[l](\cat{C})(X,Z)} \SStat(X)(C,A) \\
      \xto{\sim}\quad& \BLens_2\bigl((X,A),(Z,C)\bigr)
    \end{align*}
    where the functor in the penultimate line amounts to the product of the horizontal composition functors on $\ccopara[l](\cat{C})$ and $\SStat(X)$.
  \end{enumerate}
\end{defn}

\begin{rmk}
  When $\cat{C}$ is symmetric monoidal, $\SStat$ acquires the structure of a monoidal indexed bicategory (Definition \ref{def:mon-idx-bicat} and Theorem \ref{thm:monoidal-sstat}), and hence $\BLens_2$ becomes a monoidal bicategory (Corollary \ref{cor:mon-blens}).
\end{rmk}

\subsection{Coparameterized Bayesian updates compose optically}

So that our generalized Bayesian lenses are worthy of the name, we should also confirm that Bayesian inversions compose according to the lens pattern (`optically') also in the coparameterized setting.
Such confirmation is the subject of the present section, and so we first introduce a new ``coparameterized Bayes' rule''.

\begin{defn}
  We say that a coparameterized channel $\gamma:A\klto M\otimes B$ \textit{admits Bayesian inversion} if there exists a dually coparameterized channel $\rho_\pi:B\klto A\otimes M$ satisfying the graphical equation (with string diagrams read from bottom to top)
  \[ \scalebox{0.72}
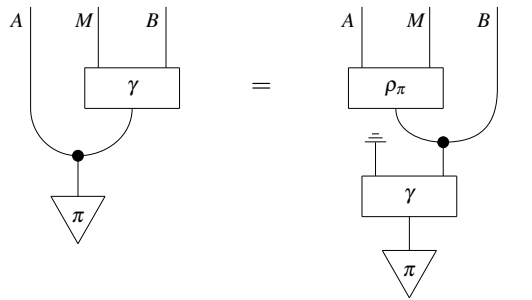
{\begin{tikzpicture}
	\begin{pgfonlayer}{nodelayer}
		\node [style=none] (0) at (0, 0.75) {};
		\node [style=none] (1) at (-0.5, 0.75) {};
		\node [style=none] (2) at (-0.5, -0.75) {};
		\node [style=none] (3) at (1.25, -0.75) {};
		\node [style=none] (4) at (3, -0.75) {};
		\node [style=none] (5) at (3, 0.75) {};
		\node [style=none] (6) at (2.5, 0.75) {};
		\node [style=none] (7) at (1.25, 0) {$\gamma$};
		\node [style=none] (8) at (-2.5, -0.75) {};
		\node [style=none] (9) at (0, 3) {};
		\node [style=none] (10) at (2.5, 3) {};
		\node [style=black copier] (11) at (-0.75, -2.5) {};
		\node [style=none] (12) at (-0.75, -4) {};
		\node [style=none] (13) at (-1.75, -4) {};
		\node [style=none] (14) at (0.25, -4) {};
		\node [style=none] (15) at (-0.75, -5.75) {};
		\node [style=none] (16) at (-0.75, -4.75) {$\pi$};
		\node [style=none] (17) at (-2.5, 3) {};
		\node [style=none] (18) at (-3, 2.5) {$A$};
		\node [style=none] (19) at (-0.5, 2.5) {$M$};
		\node [style=none] (20) at (2, 2.5) {$B$};
	\end{pgfonlayer}
	\begin{pgfonlayer}{edgelayer}
		\draw (2.center) to (1.center);
		\draw (4.center) to (5.center);
		\draw (2.center) to (4.center);
		\draw (1.center) to (5.center);
		\draw (17.center)
			 to (8.center)
			 to [in=180, out=-90] (11.center)
			 to [in=-90, out=0] (3.center);
		\draw (12.center) to (11);
		\draw (13.center) to (14.center);
		\draw (14.center) to (15.center);
		\draw (15.center) to (13.center);
		\draw (0.center) to (9.center);
		\draw (6.center) to (10.center);
	\end{pgfonlayer}
\end{tikzpicture}} \qquad = \qquad \scalebox{0.72}{\begin{tikzpicture}
	\begin{pgfonlayer}{nodelayer}
		\node [style=none] (0) at (-2.5, 0.75) {};
		\node [style=none] (1) at (-3, 0.75) {};
		\node [style=none] (2) at (-3, -0.75) {};
		\node [style=none] (3) at (-1.25, -0.75) {};
		\node [style=none] (4) at (0.5, -0.75) {};
		\node [style=none] (5) at (0.5, 0.75) {};
		\node [style=none] (6) at (0, 0.75) {};
		\node [style=none] (7) at (-1.25, 0) {$\rho_\pi$};
		\node [style=none] (8) at (-2.5, 3) {};
		\node [style=none] (9) at (0, 3) {};
		\node [style=none] (10) at (-2, -3.25) {};
		\node [style=none] (11) at (-2.5, -3.25) {};
		\node [style=none] (12) at (-2.5, -4.75) {};
		\node [style=none] (14) at (1, -4.75) {};
		\node [style=none] (15) at (1, -3.25) {};
		\node [style=none] (16) at (0.5, -3.25) {};
		\node [style=none] (17) at (-0.75, -4) {$\gamma$};
		\node [style=discarder] (18) at (-2, -2) {};
		\node [style=black copier] (19) at (0.5, -2) {};
		\node [style=none] (21) at (2.5, 0) {};
		\node [style=none] (23) at (2.5, 3) {};
		\node [style=none] (24) at (-0.75, -6) {};
		\node [style=none] (25) at (-0.75, -4.75) {};
		\node [style=none] (27) at (-1.75, -6) {};
		\node [style=none] (28) at (0.25, -6) {};
		\node [style=none] (29) at (-0.75, -7.75) {};
		\node [style=none] (30) at (-0.75, -6.75) {$\pi$};
		\node [style=none] (31) at (-3, 2.5) {$A$};
		\node [style=none] (32) at (-0.5, 2.5) {$M$};
		\node [style=none] (33) at (2, 2.5) {$B$};
	\end{pgfonlayer}
	\begin{pgfonlayer}{edgelayer}
		\draw (2.center) to (1.center);
		\draw (4.center) to (5.center);
		\draw (2.center) to (4.center);
		\draw (1.center) to (5.center);
		\draw (0.center) to (8.center);
		\draw (6.center) to (9.center);
		\draw (12.center) to (11.center);
		\draw (14.center) to (15.center);
		\draw (12.center) to (14.center);
		\draw (11.center) to (15.center);
		\draw (10.center) to (18);
		\draw (16.center) to (19);
		\draw (3.center)
			 to [in=-180, out=-90] (19.center)
			 to [in=-90, out=0, looseness=1.25] (21.center)
			 to (23.center);
		\draw (24.center) to (25.center);
		\draw (27.center) to (28.center);
		\draw (28.center) to (29.center);
		\draw (29.center) to (27.center);
	\end{pgfonlayer}
\end{tikzpicture}} \quad . \]
  In this case, we say that $\rho_\pi$ is the \textit{Bayesian inversion of} $\gamma$ \textit{with respect to} $\pi$.
\end{defn}

With this definition, we can supply the desired result that ``coparameterized Bayesian updates compose optically''.

\begin{thm} \label{thm:copara2-buco}
  Suppose $(\gamma,\gamma^\dag):(A,A)\xlensto[M]{}(B,B)$ and $(\delta,\delta^\dag):(B,B)\xlensto[N]{}(C,C)$ are coparameterized Bayesian lenses in $\BLens_2$.
  Suppose also that $\pi:I\klto A$ is a state on $A$ in the underlying category of channels $\cat{C}$, such that $\gamma^\dag_\pi$ is a Bayesian inversion of $\gamma$ with respect to $\pi$, and such that $\delta^\dag_{\gamma\pi}$ is a Bayesian inversion of $\delta$ with respect to $(\gamma\pi)^{\smallground}$; where the notation $(-)^{\smallground}$ represents discarding coparameters.
  Then $\gamma^\dag_\pi\klcirc\delta^\dag_{\gamma\pi}$ is a Bayesian inversion of $\delta\klcirc\gamma$ with respect to $\pi$.
  (Here $\klcirc$ denotes copy-composition.)
  Moreover, if $(\delta\klcirc\gamma)^\dag_\pi$ is any Bayesian inversion of $\delta\klcirc\gamma$ with respect to $\pi$, then $\gamma^\dag_\pi\klcirc\delta^\dag_{\gamma\pi}$ is $(\delta\gamma\pi)^{\smallground}$-almost-surely equal to $(\delta\klcirc\gamma)^\dag_\pi$: that is, $(\delta\klcirc\gamma)^\dag_\pi \overset{(\delta\gamma\pi)^{\smallground}}{\sim} \gamma^\dag_\pi\klcirc\delta^\dag_{\gamma\pi}$.
\end{thm}

In order to satisfy this coparameterized Bayes' rule, a Bayesian lens must be of `simple' type.

\begin{defn}
  We say that a coparameterized Bayesian lens $(c,c')$ is \textit{simple} if its domain and codomain are `diagonal' (duplicate pairs of objects) and if the coparameter of $c$ is equal to the coparameter of $c'$.
  In this case, we can write the type of $(c,c')$ as $(X,X)\xlensto[M]{}(Y,Y)$ or simply $X\xlensto[M]{}Y$.
\end{defn}

\section{Statistical games for local approximate inference} \label{sec:stat-games}

\subsection{Losses for lenses}

Statistical games are obtained by attaching to Bayesian lenses \textit{loss functions}, representing `local' quantifications of the performance of approximate inference systems.
Because this performance depends on the system's context (\textit{i.e.}, the prior $\pi:I\klto X$ and the observed data $b:B$), a loss function at its most concrete will be a function $\cat{C}(I,X)\times B\to\rr_+$.
To internalize this type in $\cat{C}$, we may recall that, when $\cat{C}$ is the category $\Cat{sfKrn}$ of s-finite kernels or the Kleisli category $\Kl(\Da_{\leq1})$ of the subdistribution monad, a density function $p_c:X\times Y\to[0,1]$ for a channel $c:X\klto Y$ corresponds to an \textit{effect} (or \textit{costate}) $X\otimes Y\klto I$.
In this way, we can see a loss function as a kind of \textit{state-dependent effect} $B\xklto{X}I$.

Loss functions will compose by sum, and so we need to ask for the effects in $\cat{C}$ to form a monoid.
Moreover, we need this monoid to be `bilinear' with respect to channels, so that $\Stat$-reindexing (\textit{cf.} Definition \ref{def:4-stat-cat}) preserves sums.
These conditions are formalized in the following definition.

\begin{defn} \label{def:bilin-eff}
  Suppose $(\cat{C},\otimes,I)$ is a copy-discard category.
  We say that $\cat{C}$ \textit{has bilinear effects} if the following conditions are satisfied:
  \begin{enumerate}[label=(\roman*)]
  \item \textit{effect monoid}: there is a natural transformation $+:\cat{C}(-,I)\times\cat{C}(=,I)\Rightarrow\cat{C}({-}\otimes{=},I)$ making $\sum_{A:\cat{C}}\cat{C}(A,I)$ into a commutative monoid with unit $0:I\klto I$;
  \item \textit{bilinearity}: $(g+g')\klcirc\bcopier\klcirc f = g\klcirc f + g'\klcirc f$ for all effects $g,g'$ and morphisms $f$ such that $(g+g')\klcirc\bcopier\klcirc f$ exists.
  \end{enumerate}
\end{defn}

A trivial example of a category with bilinear effects is supplied by any Cartesian category, such as $\Set$.
If $M$ is any monoid in $\Set$, then a less trivial example is supplied by the Kleisli category of the corresponding free module monad; bilinearity follows from the module structure.
A related non-example is $\Kl(\Da_{\leq1})$: the failure here is that the effects only form a \textit{partial} monoid\footnote{
Indeed, an \textit{effect algebra} is a kind of partial monoid \parencite[\S2]{Jacobs2015Effect}, but we do not need the extra complication here.}.
More generally, the category $\Cat{sfKrn}$ of s-finite kernels \parencite{Vakar2018S} has bilinear effects (owing to the linearity of integration), and we will assume this as our ambient $\cat{C}$ for the examples below.

Given such a category $\cat{C}$ with bilinear effects, we can lift the natural transformation $+$, and hence the bilinear effect structure, to the fibres of $\Stat_{\cat{C}}$, using the universal property of the product of categories:
\begin{align*}
  +_X &: \Stat(X)(-,I) \times \Stat(X)(=,I) \,\xlongequal{}\, \Set\bigl(\cat{C}(I,X),\cat{C}(-,I)\bigr) \times \Set\bigl(\cat{C}(I,X),\cat{C}(=,I)\bigr) \\
  &\xRightarrow{(\cdot,\cdot)} \Set\bigl(\cat{C}(I,X),\cat{C}(-,I)\times\cat{C}(=,I)\bigr) \\
  &\xRightarrow{\Set\bigl(\cat{C}(I,X),+\bigr)} \Set\bigl(\cat{C}(I,X),\cat{C}({-}\otimes{=},I)\bigr) \\
  &\xRightarrow{=} \Stat(X)({-}\otimes{=},I)
\end{align*}
Here, $(\cdot,\cdot)$ denotes the pairing operation obtained from the universal property.
In this way, each $\Stat(X)$ has bilinear effects.
Note that this lifting is (strictly) compatible with the reindexing of $\Stat$, so that $+_{(-)}$ defines an indexed natural transformation.
This means in particular that \textit{reindexing distributes over sums}: given state-dependent effects $g,g':B\xklto{Y}I$ and a channel $c:X\klto Y$, we have $(g+_Y g')_c = g_c +_X g'_c$.
We will thus generally omit the subscript from the lifted sum operation, and just write $+$.

We are now ready to construct the bicategory of statistical games.

\begin{defn} \label{def:bicat-sgame}
  Suppose $(\cat{C},\otimes,I)$ has bilinear effects, and let $\BLens_2$ denote the corresponding bicategory of (copy-composite) Bayesian lenses.
  We will write $\SGame{\cat{C}}$ to denote the following \textit{bicategory of (copy-composite) statistical games} in $\cat{C}$:
  \begin{itemize}
  \item The 0-cells are the 0-cells $(X,A)$ of $\BLens_2$;
  \item the 1-cells, called \textit{statistical games}, $(X,A)\to(Y,B)$ are pairs $(c,L^c)$ of a 1-cell $c:(X,A)\lensto(Y,B)$ in $\BLens_2$ and a \textit{loss} $L^c:B\xklto{X}I$ in $\Stat(X)(B,I)$;
  \item given 1-cells $(c,L^c),(c',L^{c'}):(X,A)\to(Y,B)$, the 2-cells $(c,L^c)\Rightarrow(c',L^{c'})$ are pairs $(\alpha,K^\alpha)$ of a 2-cell $\alpha:c\Rightarrow c'$ in $\BLens_2$ and a loss $K^\alpha:B\xklto{X}I$ such that $L^{c} = L^{c'} + K^\alpha$;
  \item the identity 2-cell on $(c,L^c)$ is $(\id_c,0)$;
  \item given 2-cells $(\alpha,K^\alpha):(c,L^c)\Rightarrow(c',L^{c'})$ and $(\alpha',K^{\alpha'}):(c',L^{c'})\Rightarrow(c'',L^{c''})$, their vertical composite is $(\alpha'\circ\alpha, K^{\alpha'}+K^{\alpha})$, where $\circ$ here denotes vertical composition in $\BLens_2$;
  \item given 1-cells $(c,L^c):(X,A)\to(Y,B)$ and $(d,L^d):(Y,B)\to(Z,C)$, their horizontal composite is $(c\lenscirc d, L^d_c + L^c\circ\overline{d}_c)$; and
    \begin{itemize}
    \item given 2-cells $(\alpha,K^\alpha):(c,L^c)\Rightarrow(c',L^{c'})$ and $(\beta,K^\beta):(d,L^d)\Rightarrow(d',L^{d'})$, their horizontal composite is $(\beta\lenscirc\alpha,K^\beta_c + K^\alpha\circ\overline{d}_c)$, where $\lenscirc$ here denotes horizontal composition in $\BLens_2$.
    \end{itemize}
  \end{itemize}
\end{defn}

\begin{thm} \label{thm:bicat-sgame-well-def}
  Definition \ref{def:bicat-sgame} generates a well-defined bicategory.
\end{thm}

The proof of this result is that $\SGame{\cat{C}}$ is obtained via a pair of bicategorical Grothendieck constructions \parencite{Bakovic2010Fibrations}: first to obtain Bayesian lenses; and then to attach the loss functions.
The proof depends on the intermediate result that our effect monoids can be `upgraded' to monoidal categories; we then use the delooping of this structure to associate (state-dependent) losses to (state-dependent) channels, after discarding the coparameters of the latter.

\begin{lemma} \label{lemma:effect-mon-cats}
  Suppose $(\cat{C},\otimes,I)$ has bilinear effects.
  Then, for each object $B$, $\cat{C}(B,I)$ has the structure of a symmetric monoidal category.
  The objects of $\cat{C}(B,I)$ are its elements, the effects.
  If $g,g'$ are two effects, then a morphism $\kappa:g\to g'$ is an effect such that $g = g' + \kappa$; the identity morphism for each effect $\id_g$ is then the constant $0$ effect.
  Likewise, the tensor of two effects is their sum, and the corresponding unit is the constant $0$.
  Precomposition by any morphism $c:A\klto B$ preserves the monoidal category structure, making the presheaf $\cat{C}(-,I)$ into a fibrewise-monoidal indexed category $\cat{C}\op\to\Cat{MonCat}$.
\end{lemma}

As already indicated, this structure lifts to the fibres of $\Stat$.

\begin{cor} \label{cor:stat-effect-mon-cats}
  For each object $X$ in a category with bilinear effects, and for each object $B$, $\Stat(X)(B,I)$ inherits the symmetric monoidal structure of $\cat{C}(B,I)$; note that morphisms of state-dependent effects are likewise state-dependent, and that tensoring (summing) state-dependent effects involves copying the parameterizing state.
  Moreover, $\Stat(-)(=,I)$ is a fibrewise-monoidal indexed category $\sum_{X:\cat{C}\op}\Stat(X)\op \to \Cat{MonCat}$.
\end{cor}

\subsection{Local inference models}

In the context of approximate inference, one often does not have a single statistical model to evaluate, but a whole family of them.
In particularly nice situations, this family is actually a subcategory $\cat{D}$ of $\cat{C}$, with the family of statistical models being all those that can be composed in $\cat{D}$.
The problem of approximate inference can then be formalized as follows.
Since both $\BLens_2$ and $\SGame{\cat{C}}$ were obtained by bicategorical Grothendieck constructions, we have a pair of 2-fibrations $\SGame{\cat{C}}\xto{\pLoss}\BLens_2\xto{\pLens}\ccopara[l](\cat{C})$.
Each of $\pLoss$, $\pLens$, and the discarding functor $(-)^{\smallground}$ can be restricted to the subcategory $\cat{D}$.
The inclusion $\cat{D}\hookrightarrow\ccopara[l](\cat{D})$ is a section of this restriction of $(-)^{\smallground}$; the assignment of inversions to channels in $\cat{D}$ then corresponds to a 2-section of the 2-fibration $\pLens$ (restricted to $\cat{D}$); and the subsequent assignment of losses is a further 2-section of $\pLoss$.
This situation is depicted in the following diagram of bicategories:
\[\begin{tikzcd}[row sep=scriptsize]
	{\SGame{\cat{D}}} & {\SGame{\cat{C}}} \\
	{\BLens_2|_{\cat{D}}} & {\BLens_2} \\
	{\ccopara[l](\cat{D})} & {\ccopara[l](\cat{C})} \\
	{\cat{D}} & {\cat{C}}
	\arrow["\pLoss", from=1-2, to=2-2]
	\arrow["\pLens", from=2-2, to=3-2]
	\arrow["\ground", from=3-2, to=4-2]
	\arrow[hook, from=4-1, to=4-2]
	\arrow[hook, from=3-1, to=3-2]
	\arrow[hook, from=2-1, to=2-2]
	\arrow[hook, from=1-1, to=1-2]
	\arrow["{\pLoss|_{\cat{D}}}", from=1-1, to=2-1]
	\arrow["{\pLens|_{\cat{D}}}", from=2-1, to=3-1]
	\arrow["{\ground|_{\cat{D}}}", from=3-1, to=4-1]
	\arrow[curve={height=-24pt}, from=2-1, to=1-1]
	\arrow[curve={height=-24pt}, from=3-1, to=2-1]
	\arrow[curve={height=-24pt}, hook', from=4-1, to=3-1]
\end{tikzcd}\]
This motivates the following definitions of \textit{inference system} and \textit{loss model}, although, for the sake of our examples, we will explicitly allow the loss-assignment to be lax.
Before giving these new definitions, we recall the notion of \textit{essential image} of a functor.

\begin{defn}[\parencite{nLab2023EssentialImage}]
  Suppose $F:\cat{C}\to\cat{D}$ is an n-functor (a possibly weak homomorphism of weak n-categories).
  The \textit{image} of $F$ is the smallest sub-n-category of $\cat{D}$ that contains $F(\alpha)$ for all k-cells $\alpha$ in $\cat{C}$, along with any $(k+1)$-cells relating images of composites and composites of images, for all $0\leq k\leq n$.
  We say that a sub-n-category $\cat{D}$ is \textit{replete} if, for any k-cells $\alpha$ in $\cat{D}$ and $\beta$ in $\cat{C}$ (with $0\leq k<n$) such that $f:\alpha\Rightarrow\beta$ is a $(k+1)$-isomorphism in $\cat{C}$, then $f$ is also a $(k+1)$-isomorphism in $\cat{D}$.
  The \textit{essential image} of $F$, denoted $\im(F)$, is then the smallest replete sub-n-category of $\cat{D}$ containing the image of $F$.
\end{defn}

\begin{defn}
  Suppose $(\cat{C},\otimes,I)$ is a copy-delete category.
  An \textit{inference system} in $\cat{C}$ is a pair $(\cat{D},\ell)$ of a subcategory $\cat{D}\hookrightarrow\cat{C}$ along with a section $\ell:\im(\iota)\to\BLens_2|_{\cat{D}}$ of $\pLens|_{\cat{D}}$, where $\im(\iota)$ is the essential image of the canonical lax inclusion $\iota:\cat{D}\hookrightarrow\ccopara[l](\cat{D})$.
\end{defn}

\begin{defn}
  Suppose $(\cat{C},\otimes,I)$ has bilinear effects and $\cat{B}$ is a subbicategory of $\BLens_2$.
  A \textit{loss model} for $\cat{B}$ is a lax section $L$ of the restriction $\pLoss|_{\cat{B}}$ of $\pLoss$ to $\cat{B}$.
  We say that $L$ is a \textit{strict} loss model if it is in fact a strict 2-functor, and a \textit{strong} loss model if it is in fact a pseudofunctor.
\end{defn}

\begin{rmk}
  We may often be interested in loss models for which $\cat{B}$ is in fact the essential image of an inference system, but we do not stipulate this requirement in the definition as it is not necessary for the following development.
\end{rmk}

Since lax functors themselves collect into categories, and using the monoidality of $+$, we have the following easy proposition that will prove useful below.

\begin{prop} \label{prop:loss-smc}
  Loss models for $\cat{B}$ constitute the objects of a symmetric monoidal category $\bigl(\Loss(\cat{B}), +, 0\bigr)$.
  The morphisms of $\Loss(\cat{B})$ are icons (identity component oplax transformations \parencite[{\S4.6}]{Johnson20202Dimensional}) between the corresponding lax functors, and they compose accordingly.
  The monoidal structure is given by sums of losses.
\end{prop}

\subsection{Examples}

Each of our examples involves taking expectations of log-densities, and so to make sense of them it first helps to understand what we mean by ``taking expectations''.

\begin{notation}[Expectations] \label{nota:expectations}
  Written as a function, a density $p$ on $X$ has the type $X\to \rr_+$; written as an effect, the type is $X\klto I$.
  Given a measure or distribution $\pi$ on $X$ (equivalently, a state $\pi:I\klto X$), we can compute the expectation of $p$ under $\pi$ as the composite $p\klcirc\pi$.
  We write the resulting quantity as $\E_\pi[p]$, or more explicitly as $\E_{x\sim\pi}\bigl[p(x)\bigr]$.
  We can think of this expectation as representing the `validity' (or truth value) of the `predicate' $p$ given the state $\pi$ \parencite{Jacobs2020Logical}.
\end{notation}

\subsubsection{Relative entropy and Bayesian inference}

For our first example, we return to the subject with which we opened this paper: the compositional structure of the relative entropy.
We begin by giving a precise definition.

\begin{defn}
  Suppose $\alpha,\beta$ are both measures on $X$, with $\alpha$ absolutely continuous with respect to $\beta$.
  Then the \textit{relative entropy} or \textit{Kullback-Leibler divergence} from $\alpha$ to $\beta$ is the quantity $D_{KL}(\alpha,\beta) := \E_\alpha \left[ \log \frac{\alpha}{\beta} \right]$, where $\frac{\alpha}{\beta}$ is the Radon-Nikodym derivative of $\alpha$ with respect to $\beta$.
\end{defn}

\begin{rmk}
  When $\alpha$ and $\beta$ admit density functions $p_\alpha$ and $p_\beta$ with respect to the same base measure $\d x$, then $D_{KL}(\alpha,\beta)$ can equally be computed as $\E_{x\sim\alpha} \bigl[ \log p_\alpha(x) - \log p_\beta(x) \bigr]$.
  It it this form that we will adopt henceforth.
\end{rmk}

\begin{prop} \label{prop:kl-loss-model}
  Let $\cat{B}$ be a subbicategory of simple lenses in $\BLens_2$, all of whose channels admit density functions with respect to a common measure and whose forward channels admit Bayesian inversion (and whose forward and backward coparameters coincide), and with only structural 2-cells.
  Then the relative entropy defines a strict loss model $\KL:\cat{B}\to\SGame{}$.
  Given a lens $(c,c'):(X,X)\lensto(Y,Y)$, $\KL$ assigns the loss function $\KL(c,c'):Y\xklto{X}I$ defined, for $\pi:I\klto X$ and $y:Y$, by the relative entropy $\KL(c,c')_\pi(y) := D_{KL}\bigl(c'_\pi(y),c^\dag_\pi(y)\bigr)$,
  where $c^\dag$ is the exact inversion of $c$.
\end{prop}

Successfully playing a relative entropy game entails minimizing the divergence from the approximate to the exact posterior.
This divergence is minimized when the two coincide, and so $\KL$ represents a form of approximate Bayesian inference.

\subsubsection{Maximum likelihood estimation}

A statistical system may be more interested in predicting observations than updating beliefs.
This is captured by the process of \textit{maximum likelihood estimation}.

\begin{defn}
  Let $(c,c'):(X,X)\lensto(Y,Y)$ be a simple lens whose forward channel $c$ admits a density function $p_c$.
  Then its \textit{log-likelihood} is the loss function defined by $\MLE(c,c')_\pi(y) := -\log p_{c^{\smallground}\klcirc\pi}(y)$.
\end{defn}

\begin{prop} \label{prop:mle-loss-model}
  Let $\cat{B}$ be a subbicategory of lenses in $\BLens_2$ all of which admit density functions with respect to a common measure, and with only structural 2-cells.
  Then the assignment $(c,c') \mapsto \MLE(c,c')$ defines a lax loss model $\MLE:\cat{B}\to\SGame{}$.
\end{prop}

Successfully playing a maximum likelihood game involves maximizing the log-likelihood that the system assigns to its observations $y:Y$.
This process amounts to choosing a channel $c$ that assigns high likelihood to likely observations, and thus encodes a valid model of the data distribution.

\subsubsection{Autoencoders via the free energy}

Many adaptive systems neither just infer nor just predict: they do both, building a model of their observations that they also invert to update their beliefs.
In machine learning, such systems are known as \textit{autoencoders}, as they `encode' (infer) and `decode' (predict), `autoassociatively' \parencite{Kramer1991Nonlinear}.
In a Bayesian context, they are known as \textit{variational autoencoders} \parencite{Kingma2013Auto}, and their loss function is the \textit{free energy} \parencite{Dayan1995Helmholtz}.

\begin{defn} \label{def:fe-loss-model}
  The \textit{free energy} loss model is the sum of the relative entropy and the likelihood loss models: $\FE := \KL + \MLE$.
  Given a simple lens $(c,c'):(X,X)\lensto(Y,Y)$ admitting Bayesian inversion and with densities, $\FE$ assigns the loss function
  \begin{align*}
    \FE(c,c')_\pi(y) &= (\KL + \MLE)(c,c')_\pi(y) \\
    &= D_{KL}\bigl(c'_\pi(y),c^\dag_\pi(y)\bigr) - \log p_{c^{\smallground}\klcirc\pi}(y)
  \end{align*}
\end{defn}

\begin{rmk}
  Beyond its autoencoding impetus, another important property of the free energy is its improved computational tractability compared to either the relative entropy or the likelihood loss.
  This property is a consequence of the following fact:
  although obtained as the sum of terms which both depend on an expensive marginalization\footnote{
  Evaluating the pushforward $c^{\smallground}\klcirc\pi$ involves marginalizing over the intermediate variable; and evaluating $c^\dag_\pi(y)$ also involves evaluating $c^{\smallground}\klcirc\pi$.},
  the free energy itself does not.
  This can be seen by expanding the definitions of the relative entropy and of $c^\dag_\pi$ and rearranging terms:
  \begin{align*}
    \FE(c,c')_\pi(y)
    &= D_{KL}\bigl(c'_\pi(y),c^\dag_\pi(y)\bigr) - \log p_{c^{\smallground}\klcirc\pi}(y) \\
    &= \E_{(x,m)\sim c'_\pi(y)} \bigl[ \log p_{c'_\pi}(x,m|y) - \log p_{c^\dag_\pi}(x,m|y) \bigr] - \log p_{c^{\smallground}\klcirc\pi}(y) \\
    &= \E_{(x,m)\sim c'_\pi(y)} \bigl[ \log p_{c'_\pi}(x,m|y) - \log p_{c^\dag_\pi}(x,m|y) - \log p_{c^{\smallground}\klcirc\pi}(y) \bigr] \\
    &= \E_{(x,m)\sim c'_\pi(y)} \bigl[ \log p_{c'_\pi}(x,m|y) - \log \frac{p_c(m,y|x) p_\pi(x)}{p_{c^{\smallground}\klcirc\pi}(y)} - \log p_{c^{\smallground}\klcirc\pi}(y) \bigr] \\
    &= \E_{(x,m)\sim c'_\pi(y)} \bigl[ \log p_{c'_\pi}(x,m|y) - \log p_c(m,y|x) - \log p_\pi(x) \bigr] \\
    &= D_{KL}\bigl(c'_\pi(y), \pi\otimes\mathbf{1}\bigr) - \E_{(x,m)\sim c'_\pi(y)} \bigl[ \log p_c(m,y|x) \bigr]
  \end{align*}
  Here, $\mathbf{1}$ denotes the measure with density $1$ everywhere.
  Note that when the coparameter is trivial, $\FE(c,c')_\pi(y)$ reduces to
  \[ D_{KL}\bigl(c'_\pi(y), \pi\bigr) - \E_{x\sim c'_\pi(y)} \bigl[ \log p_c(y|x) \bigr] \; . \]
\end{rmk}

\begin{rmk}
  The name \textit{free energy} is due to an analogy with the Helmholtz free energy in thermodynamics, as we can write it as the difference between an (expected) energy and an entropy term:
  \begin{align*}
    \FE(c,c')_\pi(y)
    &= \E_{(x,m) \sim c'_\pi(y)} \bigl[ - \log p_c(m,y|x) - \log p_\pi (x) \bigr]
       - S_{X\otimes M} \bigl[ c'_\pi(y) \bigr] \\
    &= \E_{(x,m) \sim c'_\pi(y)} \left[ E_{(c,\pi)}(x,m,y) \right] - S_{X\otimes M} \bigl[ c'_\pi(y) \bigr]
  \; = U - TS
  \end{align*}
  where we call $E_{(c,\pi)} : X\otimes M \otimes Y \xklto{X} I$ the \textit{energy}, and where $S_{X\otimes M} : I\xklto{X\otimes M}I$ is the Shannon entropy.
  The last equality makes the thermodynamic analogy: \(U\) here is the \textit{internal energy} of the system; \(T = 1\) is the \textit{temperature}; and \(S\) is again the entropy.
\end{rmk}

\subsubsection{The Laplace approximation}

Although optimizing the free energy does not necessitate access to exact inversions, it is still necessary to compute an expectation under the approximate inversion, and unless one chooses wisely\footnote{
In machine learning, optimizing variational autoencoders using stochastic gradient descent typically involves a ``reparameterization trick'' \parencite[\S2.5]{Kingma2017Variational}.}, this might also be difficult.
One such wise choice established in the computational neuroscience literature is the Laplace approximation \parencite{Friston2007Variational}, in which one assumes Gaussian channels and posteriors with small variance.
Under these conditions, the expectations can be approximated away.

\begin{defn}
  We will say that a channel $c$ is \textit{Gaussian} if $c(x)$ is a Gaussian measure for every $x$ in its domain.
  We will denote the mean and variance of $c(x)$ by $\mu_c(x)$ and $\Sigma_c(x)$ respectively.
\end{defn}

\begin{prop}[Laplace approximation] \label{lemma:6-laplace-approx}
  Let the ambient category of channels $\cat{C}$ be restricted to that generated by Gaussian channels between finite-dimensional Cartesian spaces, and let $\cat{B}$ denote the corresponding restriction of $\BLens_2$.
  Suppose $(\gamma,\rho):(X,X)\lensto(Y,Y)$ is such a lens, for which, for all $y:Y$ and Gaussian priors $\pi:I\klto X$, the eigenvalues of $\Sigma_{\rho_\pi}(y)$ are small.
  Then the free energy $\FE(\gamma,\rho)_\pi(y)$ can be approximated by the \textit{Laplacian free energy}
  \begin{align} \label{eq:6-laplace-energy}
    \FE(\gamma,\rho)_\pi(y) & \approx \LFE(\gamma,\rho)_\pi(y) \\
    & := E_{(\gamma,\pi)}\bigl(\mu_{\rho_\pi}(y), y\bigr) - S_{X\otimes M} \bigl[ \rho_\pi(y) \bigr] \\
    & = -\log p_\gamma(\mu_{\rho_\pi}(y),y) - \log p_\pi(\mu_{\rho_\pi}(y)|_X) - S_{X\otimes M} \bigl[ \rho_\pi(y) \bigr] \nonumber
  \end{align}
  where we have written the argument of the density $p_\gamma$ in `function' style; where $(-)_X$ denotes the projection onto $X$; and where \(S_{X\otimes M}[\rho_\pi(y)] = \E_{(x,m) \sim \rho_\pi(y)} [ -\log p_{\rho_\pi}(x,m|y) ]\) is the Shannon entropy of \(\rho_\pi(y)\).
  The approximation is valid when \(\Sigma_{\rho_\pi}\) satisfies
  \begin{equation} \label{eq:6-laplace-sigma-rho-pi}
    \Sigma_{\rho_\pi} (y) = \left(\partial_{(x,m)}^2 E_{(\gamma,\pi)}\right)\left( \mu_{\rho_\pi}(y), y\right)^{-1} \, .
  \end{equation}
  We call $E_{(\gamma,\pi)}$ the \textit{Laplacian energy}.
\end{prop}

\begin{rmk}
  The usual form of the Laplace model in the literature omits the coparameters.
  It is of course easy to recover the non-coparameterized form by taking $M=1$.
\end{rmk}

\begin{prop} \label{prop:lfe-loss-model}
  Let $\cat{B}$ be a subbicategory of $\BLens_2$ of Gaussian lenses whose backward channels have small variance.
  Then $\LFE$ defines a lax loss model $\cat{B}\to\SGame{}$.
\end{prop}

Effectively, this proposition says that, under the stated conditions, the free energy and the Laplacian free energy coincide.
Consequently, successfully playing a Laplacian free energy game has the same autoencoding effect.

\section{Future work}

This paper only scratches the surface of the structure of statistical games.
One avenue for further investigation is the link between this structure and the similar structure of diegetic open (economic) games \parencite{Capucci2023Diegetic}, a recent reformulation of compositional game theory \parencite{Ghani2018Compositional}.
Notably, the composition rule for loss functions appears closely related to the Bellman equation, suggesting that algorithms for approximate inference (such as expectation-maximization) and reinforcement learning (such as backward induction) are more than superficially similar.

Another avenue for further investigation concerns mathematical neatness.
First, we seek an abstract characterization of copy-composition and $\ccopara$; it has been suggested to us\footnote{This suggestion is due to Owen Lynch.} that the computation by compilers of ``static single-assignment form'' \parencite{Kelsey1995Correspondence} by compilers may have a similar structure.
Second, the explicit constraint defining simple coparameterized Bayesian lenses is inelegant; we expect that using dependent optics \parencite{Braithwaite2021Fibre,Vertechi2023Dependent,Capucci2022Seeing} may help to encode this constraint in the type signature, at the cost of higher-powered mathematical machinery.
Finally, we seek further examples of loss models, and more abstract (and hopefully universal) characterizations of those we already have; for example, it is known that the Shannon entropy has a topological origin \parencite{Bradley2021Entropy} via a ``nonlinear derivation'' \parencite{Leinster2021Entropy}, and we expect that we can follow this connection further.

\section{References}
\printbibliography[heading=none]

\appendix
\section{State-dependent channels} \label{sec:stat}

In this section, we review the indexed category $\Stat:\cat{C}\op\to\Cat{Cat}$ of state-dependent channels in $\cat{C}$, from which Bayesian lenses are obtained.
We can think of $\Stat$ as a decategorified, non-coparameterized, version of $\SStat$, in which the hom-sets $\Stat(X)(A,B)$ of each fibre are given by $\Set\bigl(\cat{C}(I,X),\cat{C}(A,B)\bigr)$.
Reindexing is again by pre-composition, although simplified as there are now no coparameters to discard.

\begin{defn} \label{def:4-stat-cat}
  Let $(\cat{C}, \otimes, I)$ be a monoidal category.
  Define the $\cat{C}$\textit{-state-indexed} category $\Stat: \cat{C}\op \to \Cat{Cat}$ as follows.
  \begin{align}
    \Stat \;\; : \;\; \cat{C}\op \; & \to \; \Cat{Cat} \nonumber \\
    X & \mapsto \Stat(X) := \quad \begin{pmatrix*}[l]
      & \Stat(X)_0 & := \quad \;\;\; \cat{C}_0 \\
      & \Stat(X)(A, B) & := \quad \;\;\; \Set\bigl(\cat{C}(I, X), \cat{C}(A, B)\bigr) \\
      \id_A \: : & \Stat(X)(A, A) & := \quad
      \left\{ \begin{aligned}
        \id_A : & \; \cat{C}(I, X)     \to     \cat{C}(A, A) \\
        & \quad\;\;\: \rho \quad \mapsto \quad \id_A
      \end{aligned} \right. \label{eq:4-stat} \\
    \end{pmatrix*} \\ \nonumber \\
    f : \cat{C}(Y, X) & \mapsto \begin{pmatrix*}[c]
      \Stat(f) \; : & \Stat(X) & \to & \Stat(Y) \vspace*{0.5em} \\
      & \Stat(X)_0 & = & \Stat(Y)_0 \vspace*{0.5em} \\
      & \Set(\cat{C}(I, X), \cat{C}(A, B)) & \to & \Set\bigl(\cat{C}(I, Y), \cat{C}(A, B)\bigr) \vspace*{0.125em} \\
      & \alpha & \mapsto & f^\ast \alpha : \big( \, \sigma : \cat{C}(I, Y) \, \big) \mapsto \big( \, \alpha(f \klcirc \sigma) : \cat{C}(A, B) \, \big)
    \end{pmatrix*} \nonumber
  \end{align}
  Composition in each fibre $\Stat(X)$ is as in $\cat{C}$.
  Explicitly, indicating morphisms $\cat{C}(I, X) \to \cat{C}(A, B)$ in $\Stat(X)$ by $A \xklto{X} B$, and given $\alpha : A \xklto{X} B$ and $\beta : B \xklto{X} C$, their composite is $\beta \circ \alpha : A \xklto{X} C := \rho \mapsto \beta(\rho) \klcirc \alpha(\rho)$, where here we indicate composition in $\cat{C}$ by $\klcirc$ and composition in the fibres $\Stat(X)$ by $\circ$.
  Given $f : Y \klto X$ in $\cat{C}$, the induced functor $\Stat(f) : \Stat(X) \to \Stat(Y)$ acts by pre-composition.
\end{defn}

The category of non-coparameterized Bayesian lenses is then obtained as the (1-categorical) Grothendieck construction of the pointwise opposite of $\Stat$, following \textcite{Spivak2019Generalized}.

\section{Monoidal statistical games} \label{sec:mon-stat-games}

In this section, we exhibit the monoidal structures on copy-composite Bayesian lenses, statistical games, and loss models, as well as demonstrating that each of our loss models is accordingly monoidal.
We begin by expressing the monoidal structure on $\ccopara(\cat{C})$.

\begin{prop} \label{prop:mon-ccopara}
  If the monoidal structure on $\cat{C}$ is symmetric, then $\ccopara(\cat{C})$ inherits a monoidal structure $(\otimes,I)$, with the same unit object $I$ as in $\cat{C}$.
  On 1-cells $f:A\xto[M]{}B$ and $f':A'\xto[M']{}B'$, the tensor $f\otimes f':A\otimes A'\xto[M\otimes M']{}B\otimes B'$ is defined by
  \[ \scalebox{0.75}
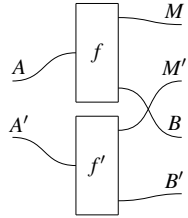
{\begin{tikzpicture}
	\begin{pgfonlayer}{nodelayer}
		\node [style=none] (0) at (0.75, 3.25) {};
		\node [style=none] (1) at (0.75, 3.75) {};
		\node [style=none] (2) at (-0.75, 3.75) {};
		\node [style=none] (3) at (-0.75, 2) {};
		\node [style=none] (4) at (-0.75, 0.25) {};
		\node [style=none] (5) at (0.75, 0.25) {};
		\node [style=none] (6) at (0.75, 0.75) {};
		\node [style=none] (7) at (0, 2) {$f$};
		\node [style=none] (8) at (-3, 1) {};
		\node [style=none] (9) at (-2.75, 1.5) {$A$};
		\node [style=none] (10) at (0.75, -0.75) {};
		\node [style=none] (11) at (0.75, -0.25) {};
		\node [style=none] (12) at (-0.75, -0.25) {};
		\node [style=none] (13) at (-0.75, -2) {};
		\node [style=none] (14) at (-0.75, -3.75) {};
		\node [style=none] (15) at (0.75, -3.75) {};
		\node [style=none] (16) at (0.75, -3.25) {};
		\node [style=none] (17) at (0, -2) {$f'$};
		\node [style=none] (18) at (-3, -1) {};
		\node [style=none] (19) at (-2.75, -0.5) {$A'$};
		\node [style=none] (20) at (3, 1) {};
		\node [style=none] (21) at (3, 3) {};
		\node [style=none] (22) at (3, -1) {};
		\node [style=none] (23) at (3, -3) {};
		\node [style=none] (24) at (2.75, 3.5) {$M$};
		\node [style=none] (25) at (2.75, 1.5) {$M'$};
		\node [style=none] (26) at (2.75, -0.5) {$B$};
		\node [style=none] (27) at (2.75, -2.5) {$B'$};
	\end{pgfonlayer}
	\begin{pgfonlayer}{edgelayer}
		\draw [in=-180, out=0, looseness=1.25] (8.center) to (3.center);
		\draw (2.center) to (1.center);
		\draw (4.center) to (5.center);
		\draw (2.center) to (4.center);
		\draw (1.center) to (5.center);
		\draw [in=-180, out=0, looseness=1.25] (18.center) to (13.center);
		\draw (12.center) to (11.center);
		\draw (14.center) to (15.center);
		\draw (12.center) to (14.center);
		\draw (11.center) to (15.center);
		\draw [in=180, out=0, looseness=1.25] (0.center) to (21.center);
		\draw [in=-180, out=0, looseness=1.25] (10.center) to (20.center);
		\draw [in=-180, out=0, looseness=1.25] (6.center) to (22.center);
		\draw [in=180, out=0, looseness=1.25] (16.center) to (23.center);
	\end{pgfonlayer}
\end{tikzpicture}} \quad . \]
  On 2-cells $\varphi:f\Rightarrow g$ and $\varphi':f'\Rightarrow g'$, the tensor $\varphi\otimes\varphi':(f\otimes f')\Rightarrow(g\otimes g')$ is given by the string diagram
  \[ \scalebox{0.75}
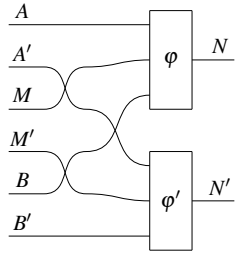
{\begin{tikzpicture}
	\begin{pgfonlayer}{nodelayer}
		\node [style=none] (0) at (1, 2.5) {};
		\node [style=none] (1) at (2.5, 2.5) {};
		\node [style=none] (2) at (1.75, 2.5) {$\varphi$};
		\node [style=none] (3) at (1, 3.75) {};
		\node [style=none] (4) at (1, 4.25) {};
		\node [style=none] (5) at (1, 0.75) {};
		\node [style=none] (6) at (2.5, 4.25) {};
		\node [style=none] (7) at (2.5, 0.75) {};
		\node [style=none] (8) at (1, 1.25) {};
		\node [style=none] (9) at (1, -2.5) {};
		\node [style=none] (10) at (2.5, -2.5) {};
		\node [style=none] (11) at (1.75, -2.5) {$\varphi'$};
		\node [style=none] (12) at (1, -1.25) {};
		\node [style=none] (13) at (1, -0.75) {};
		\node [style=none] (14) at (1, -4.25) {};
		\node [style=none] (15) at (2.5, -0.75) {};
		\node [style=none] (16) at (2.5, -4.25) {};
		\node [style=none] (17) at (1, -3.75) {};
		\node [style=none] (18) at (-1.25, -0.75) {};
		\node [style=none] (19) at (-1.25, 0.75) {};
		\node [style=none] (20) at (-1.25, 2.25) {};
		\node [style=none] (21) at (-2.75, 2.25) {};
		\node [style=none] (22) at (-2.75, 0.75) {};
		\node [style=none] (23) at (-2.75, -0.75) {};
		\node [style=none] (24) at (-2.75, -2.25) {};
		\node [style=none] (25) at (-1.25, -2.25) {};
		\node [style=none] (26) at (-4, 2.25) {};
		\node [style=none] (27) at (-4, 0.75) {};
		\node [style=none] (28) at (-4, 3.75) {};
		\node [style=none] (29) at (-4, -0.75) {};
		\node [style=none] (30) at (-4, -2.25) {};
		\node [style=none] (31) at (-4, -3.75) {};
		\node [style=none] (32) at (4, 2.5) {};
		\node [style=none] (33) at (4, -2.5) {};
		\node [style=none] (34) at (-3.5, 4.25) {$A$};
		\node [style=none] (35) at (-3.5, 2.75) {$A'$};
		\node [style=none] (36) at (-3.5, 1.25) {$M$};
		\node [style=none] (37) at (-3.5, -0.25) {$M'$};
		\node [style=none] (38) at (-3.5, -1.75) {$B$};
		\node [style=none] (39) at (-3.5, -3.25) {$B'$};
		\node [style=none] (40) at (3.5, -2) {$N'$};
		\node [style=none] (41) at (3.5, 3) {$N$};
	\end{pgfonlayer}
	\begin{pgfonlayer}{edgelayer}
		\draw (5.center) to (4.center);
		\draw (6.center) to (7.center);
		\draw (5.center) to (7.center);
		\draw (4.center) to (6.center);
		\draw (14.center) to (13.center);
		\draw (15.center) to (16.center);
		\draw (14.center) to (16.center);
		\draw (13.center) to (15.center);
		\draw (26.center)
			 to (21.center)
			 to [in=-180, out=0, looseness=1.25] (19.center)
			 to [in=-180, out=0, looseness=1.25] (12.center);
		\draw (30.center)
			 to (24.center)
			 to [in=180, out=0, looseness=1.25] (18.center)
			 to [in=-180, out=0, looseness=1.25] (8.center);
		\draw (27.center)
			 to (22.center)
			 to [in=-180, out=0, looseness=1.25] (20.center)
			 to [in=180, out=0, looseness=1.25] (0.center);
		\draw (29.center)
			 to (23.center)
			 to [in=-180, out=0, looseness=1.25] (25.center)
			 to [in=180, out=0, looseness=1.25] (9.center);
		\draw (28.center) to (3.center);
		\draw (31.center) to (17.center);
		\draw (10.center) to (33.center);
		\draw (1.center) to (32.center);
	\end{pgfonlayer}
\end{tikzpicture}} \quad . \]
\end{prop}

Next, we define the notion of monoidal indexed bicategory.

\begin{defn} \label{def:mon-idx-bicat}
  Suppose $(\cat{B},\otimes,I)$ is a monoidal bicategory.
  We will say that $F:\cat{B}\coop\to\Cat{Bicat}$ is a \textit{monoidal indexed bicategory} when it is equipped with the structure of a weak monoid object in the 3-category of indexed bicategories, indexed pseudofunctors, indexed pseudonatural transformations, and indexed modifications.

  More explicitly, we will take $F$ to be a monoidal indexed bicategory when it is equipped with
  \begin{enumerate}[label=(\roman*)]
  \item an indexed pseudofunctor $\mu:F(-)\times F(=)\to F(-\otimes=)$ called the \textit{multiplication}, \textit{i.e.},
    \begin{enumerate}[label=(\alph*)]
    \item a family of pseudofunctors $\mu_{X,Y}: FX\times FY \to F(X\otimes Y)$, along with
    \item for any 1-cells $f:X\to X'$ and $g:Y\to Y'$ in $\cat{B}$, a pseudonatural isomorphism $\mu_{f,g}:\mu_{X',Y'}\circ(Ff\times Fg)\Rightarrow F(f\otimes g)\circ\mu_{X,Y}$;
    \end{enumerate}
  \item a pseudofunctor $\eta:\Cat{1}\to FI$ called the \textit{unit};
  \end{enumerate}
  as well as three indexed pseudonatural isomorphisms --- an associator, a left unitor, and a right unitor --- which satisfy weak analogues of the coherence conditions for a monoidal indexed category \parencite[\S3.2]{Moeller2018Monoidal}, up to invertible indexed modifications.
\end{defn}

Using this notion, we can establish that $\SStat$ is monoidal.
(So as to demonstrate the structure, we do not omit the proof sketch.)

\begin{thm} \label{thm:monoidal-sstat}
  $\SStat$ is a monoidal indexed bicategory.
  \begin{proof}[Proof sketch]
    The multiplication $\mu$ is given first by the family of pseudofunctors $\mu_{X,Y}:\SStat(X)\times\SStat(Y)\to\SStat(X\otimes Y)$ which are defined on objects simply by tensor
    \[ \mu_{X,Y}(A,B) = A\otimes B \]
    since the objects do not vary between the fibres of $\SStat$, and on hom categories by the functors
    \begin{align*}
      & \SStat(X)(A,B) \times \SStat(Y)(A',B') \\
      &= \Cat{Cat}\bigl(\disc\cat{C}(I,X),\ccopara[r](\cat{C})(A,B)\bigr) \times \Cat{Cat}\bigl(\disc\cat{C}(I,Y),\ccopara[r](\cat{C})(A',B')\bigr) \\
      &\cong \Cat{Cat}\bigl(\disc\cat{C}(I,X)\times\disc\cat{C}(I,Y),\ccopara[r](\cat{C})(A,B)\times\ccopara[r](\cat{C})(A',B')\bigr) \\
      &\xto{\Cat{Cat}(\disc\cat{C}(I,\proj_X)\times\disc\cat{C}(I,\proj_Y), \otimes)} \Cat{Cat}\bigl(\disc\cat{C}(I,X\otimes Y)^2, \ccopara[r](\cat{C})(A\otimes A',B\otimes B')\bigr) \\
      &\xto{\Cat{Cat}(\bcopier,\id)} \Cat{Cat}(\disc\cat{C}(I,X\otimes Y),\ccopara[r](C)(A\otimes A', B\otimes B') \\
      &= \SStat(X\otimes Y)(A\otimes A', B\otimes B') \; .
    \end{align*}
    where $\Cat{Cat}\left(\bcopier,\id\right)$ indicates pre-composition with the universal (Cartesian) copying functor.
    For all $f:X\to X'$ and $g:Y\to Y'$ in $\ccopara[l](\cat{C})$, the pseudonatural isomorphisms
    \[ \mu_{f,g}:\mu_{X',Y'}\circ\bigl(\SStat(f)\times\SStat(g)\bigr) \Rightarrow \SStat(f\otimes g)\circ\mu_{X,Y} \]
    are obtained from the universal property of the product $\times$ of categories.
    The unit $\eta:\Cat{1}\to\SStat(I)$ is the pseudofunctor mapping the unique object of $\Cat{1}$ to the monoidal unit $I$.
    Associativity and unitality of this monoidal structure follow from the functoriality of the construction, given the monoidal structures on $\cat{C}$ and $\Cat{Cat}$.
  \end{proof}
\end{thm}

Just as the monoidal Grothendieck construction induces a monoidal structure on categories of lenses for monoidal pseudofunctors \parencite{Moeller2018Monoidal}, we obtain a monoidal structure on the bicategory of copy-composite bayesian lenses.

\begin{cor} \label{cor:mon-blens}
  The bicategory of copy-composite Bayesian lenses $\BLens_2$ is a monoidal bicategory.
  The monoidal unit is the object $(I,I)$.
  The tensor $\otimes$ is given on 0-cells by $(X,A)\otimes(X',A') := (X\otimes X', A\otimes A')$, and on hom-categories by
  \begin{align*}
    & \BLens_2\bigl((X,A),(Y,B)\bigr) \times \BLens_2\bigl((X,A),(Y,B)\bigr) \\
    &= \ccopara[l](\cat{C})(X,Y)\times\SStat(X)(B,A) \times \ccopara[l](\cat{C})(X',Y')\times\SStat(X')(B',A') \\
    &\xto{\sim} \ccopara[l](\cat{C})(X,Y)\times\ccopara[l](\cat{C})(X',Y') \times \SStat(X)(B,A)\times\SStat(X')(B',A')\\
    &\xto{\otimes\, \times\, \mu_{X,X'}\op} \ccopara[l](\cat{C})(X\otimes X', Y\otimes Y') \times \SStat(X\otimes X')(B\otimes B',A\otimes A') \\
    &= \BLens_2\bigl((X,A)\otimes(X',A'), (Y,B)\otimes(Y',B')\bigr) \; .
  \end{align*}
\end{cor}

And similarly, we obtain a monoidal structure on statistical games.

\begin{prop}
  The bicategory of copy-composite statistical games $\SGame{}$ is a monoidal bicategory.
  The monoidal unit is the object $(I,I)$.
  The tensor $\otimes$ is given on 0-cells as for the tensor of Bayesian lenses, and on hom-categories by
  \begin{align*}
    & \SGame{}\bigl((X,A),(Y,B)\bigr) \times \SGame{}\bigl((X',A'),(Y',B')\bigr) \\
    &= \BLens_2\bigl((X,A),(Y,B)\bigr) \times \Stat(X)(B,I) \\
    &\qquad \times \BLens_2\bigl((X',A'),(Y',B')\bigr) \times \Stat(X')(B',I) \\
    &\xto{\sim} \BLens_2\bigl((X,A),(Y,B)\bigr) \times \BLens_2\bigl((X',A'),(Y',B')\bigr) \\
    &\qquad \times \Stat(X)(B,I) \times \Stat(X')(B',I) \\
    &\xto{\otimes\,\times\,\mu_{X,X'}} \BLens_2\bigl((X,A)\otimes(X',A'),(Y,B)\otimes(Y',B')\bigr) \times \Stat(X\otimes X')(B\otimes B', I\otimes I) \\
    &\xto{\sim} \SGame{}\bigl((X,A)\otimes(X',A'),(Y,B)\otimes(Y',B')\bigr)
  \end{align*}
  where here $\mu$ indicates the multiplication of the monoidal structure on $\Stat$ \parencite[{Prop. 4.3.5}]{Smithe2022Mathematical}.
\end{prop}

We give natural definitions of monoidal inference system and monoidal loss model, which we elaborate below.

\begin{defn}
  A \textit{(lax) monoidal inference system} is an inference system $(\cat{D},\ell)$ for which $\ell$ is a lax monoidal pseudofunctor.
  A \textit{(lax) monoidal loss model} is a loss model $L$ which is a lax monoidal lax functor.
\end{defn}

\begin{rmk} \label{rmk:mon-loss-data}
  We say `lax' whenever a morphism (of any structure) only weakly preserves a monoidal operation such as composition of any order; this includes as a special case lax monoidal functors (since a monoidal category is a one-object bicategory).
  In this respect, we differ from \parencite[{\S2.2}]{Moeller2018Monoidal}, who use `weak' in the latter case; we prefer to maintain consistency.
  Following \parencite[{Def. 4.2.1}]{Johnson20202Dimensional}, we will continue to say \textit{lax} when the witness to laxness maps composites of images to images of composites (and \textit{oplax} when the witness maps images of composites to composites of images).

  These conventions mean that a loss model $L:\cat{B}\to\SGame{}$ is lax monoidal when it is equipped with strong transformations
  \[\begin{tikzcd}
    {\cat{B}\times\cat{B}} & {\SGame{}\times\SGame{}} \\
    {\cat{B}} & {\SGame{}}
    \arrow["{\otimes_{\cat{B}}}"', from=1-1, to=2-1]
    \arrow["{\otimes_{\Cat{G}}}", from=1-2, to=2-2]
    \arrow["{L\times L}", from=1-1, to=1-2]
    \arrow["L"', from=2-1, to=2-2]
    \arrow["\lambda"', shorten <=12pt, shorten >=12pt, Rightarrow, from=1-2, to=2-1]
  \end{tikzcd}
  \quad\text{and}\qquad
  \begin{tikzcd}
    {\Cat{1}} \\
    {\cat{B}} & {\SGame{}}
    \arrow[""{name=0, anchor=center, inner sep=0}, "{(I,I)}", curve={height=-20pt}, from=1-1, to=2-2]
    \arrow["{(I,I)}"', from=1-1, to=2-1]
    \arrow["L"', from=2-1, to=2-2]
    \arrow["{\lambda_0}"', shorten <=8pt, shorten >=6pt, Rightarrow, from=0, to=2-1]
  \end{tikzcd}\]
  where $\otimes_{\cat{B}}$ and $\otimes_{\Cat{G}}$ denote the monoidal products on $\cat{B}\hookrightarrow\BLens_2$ and $\SGame{}$ respectively, and when $\lambda$ and $\lambda_0$ are themselves equipped with invertible modifications satisfying coherence axioms, as in \textcite[{\S2.2}]{Moeller2018Monoidal}.

  Note that, because $L$ must be a (lax) section of the 2-fibration $\pLoss|_{\cat{B}}:\SGame{}|_{\cat{B}}\to\cat{B}$, the unitor $\lambda_0$ is forced to be trivial, picking out the identity on the monoidal unit $(I,I)$.
  Likewise, the laxator $\lambda:L(-)\otimes L(=)\Rightarrow L({-}\otimes{=})$ must have 1-cell components which are identities:
  \[ L(X,A)\otimes L(X',A) = (X,A)\otimes(X',A') = (X\otimes X', A\otimes A') = L\bigl((X,A)\otimes L(X',A)\bigr) \]
  The interesting structure is therefore entirely in the 2-cells.
  We follow the convention of \parencite[{Def. 4.2.1}]{Johnson20202Dimensional} that a strong transformation is a lax transformation with invertible 2-cell components.
  Supposing that $(c,c'):(X,A)\lensto(Y,B)$ and $(d,d'):(X',A')\lensto(Y',B')$ are 1-cells in $\cat{B}$, the corresponding 2-cell component of $\lambda$ has the form $\lambda_{c,d}:L\bigl((c,c')\otimes(d,d')\bigr)\Rightarrow L(c,c')\otimes L(d,d')$, hence filling the following square in $\SGame{}$:
  \[\begin{tikzcd}[row sep=scriptsize]
	{(X,A)\otimes(X',A')} && {(Y,B)\otimes(Y',B')} \\
	\\
	{(X,A)\otimes(X',A')} && {(Y,B)\otimes(Y',B')}
	\arrow["{L(c,c')\otimes L(d,d')}", from=1-1, to=1-3]
	\arrow[Rightarrow, no head, from=1-3, to=3-3]
	\arrow[Rightarrow, no head, from=1-1, to=3-1]
	\arrow["{L((c,c')\otimes(d,d'))}"', from=3-1, to=3-3]
	\arrow["{\lambda_{c,d}}", shorten <=18pt, shorten >=18pt, Rightarrow, from=3-1, to=1-3]
  \end{tikzcd}\]
  Intuitively, these 2-cells witness the failure of the tensor $L(c,c')\otimes L(d,d')$ of the parts to account for correlations that may be evident to the ``whole system'' $L\bigl((c,c')\otimes(d,d')\bigr)$.
\end{rmk}

Just as we have monoidal lax functors, we can have monoidal lax transformations; again, see \parencite[{\S2.2}]{Moeller2018Monoidal}.

\begin{prop}
  Monoidal loss models and monoidal icons form a subcategory $\MonLoss(\cat{B})$ of $\Loss(\cat{B})$, and the symmetric monoidal structure $(+,0)$ on the latter restricts to the former.
\end{prop}

\subsection{Examples}

In this section, we present the monoidal structure on the loss models considered above.
Because loss models $L$ are (lax) sections, following Remark \ref{rmk:mon-loss-data}, this monoidal structure is given in each case by a lax natural family of 2-cells $\lambda_{c,d}:L\bigl((c,c')\otimes(d,d')\bigr)\Rightarrow L(c,c')\otimes L(d,d')$, for each pair of lenses $(c,c'):(X,A)\lensto(Y,B)$ and $(d,d'):(X',A')\lensto(Y',B')$.
Such a 2-cell $\lambda_{c,d}$ is itself given by a loss function of type $B\otimes B'\xklto{X\otimes X'}I$ satisfying the equation $L\bigl((c,c')\otimes(d,d')\bigr) = L(c,c')\otimes L(d,d') + \lambda_{c,d}$.
Following \parencite[{Eq. 4.2.3}]{Johnson20202Dimensional}, lax naturality requires that $\lambda$ satisfy the following equation of 2-cells, where $K$ denotes the laxator (with respect to horizontal composition $\diamond$) with components $K(e,c):Le\diamond Lc \Rightarrow L(e\lenscirc c)$:
\begin{gather*}
  \begin{tikzcd}[ampersand replacement=\&,row sep=scriptsize]
    \& {(Y,B)\otimes(Y',B')} \\
    \\
    {(X,A)\otimes(X',A')} \&\& {(Z,C)\otimes(Z',C')}
    \arrow["{L(c\otimes d)}", curve={height=-12pt}, from=3-1, to=1-2]
    \arrow["{L(e\otimes f)}", curve={height=-12pt}, from=1-2, to=3-3]
    \arrow[""{name=0, anchor=center, inner sep=0}, "{L\bigl((e\lenscirc c)\otimes(f\lenscirc d)\bigr)}"{description}, from=3-1, to=3-3]
    \arrow[""{name=1, anchor=center, inner sep=0}, "{L(e\lenscirc c)\otimes L(f\lenscirc d)}"', curve={height=60pt}, from=3-1, to=3-3]
    \arrow["{K(e\otimes f,c\otimes d)\,}"', shorten >=10pt, Rightarrow, from=1-2, to=0]
    \arrow["{\lambda(e\lenscirc c,f\lenscirc d)\,}"', shorten <=9pt, shorten >=6pt, Rightarrow, from=0, to=1]
  \end{tikzcd}
  \\ = \\
  \begin{tikzcd}[ampersand replacement=\&,row sep=scriptsize]
    \&\& {(Y,B)\otimes(Y',B')} \\
    \\
    {(X,A)\otimes(X',A')} \&\& {(Y,B)\otimes(Y',B')} \&\& {(Z,C)\otimes(Z',C')}
    \arrow[""{name=0, anchor=center, inner sep=0}, "{L(c\otimes d)}"{pos=0.667}, curve={height=-12pt}, from=3-1, to=1-3]
    \arrow[""{name=1, anchor=center, inner sep=0}, "{L(e\otimes f)}"{pos=0.333}, curve={height=-12pt}, from=1-3, to=3-5]
    \arrow[""{name=2, anchor=center, inner sep=0}, "{L(e\lenscirc c)\otimes L(f\lenscirc d)}"', curve={height=60pt}, from=3-1, to=3-5]
    \arrow[Rightarrow, no head, from=1-3, to=3-3]
    \arrow[""{name=3, anchor=center, inner sep=0}, "{Lc\otimes Ld}"{pos=0.6}, from=3-1, to=3-3]
    \arrow[""{name=4, anchor=center, inner sep=0}, "{Le\otimes Lf}"{pos=0.4}, from=3-3, to=3-5]
    \arrow["{\lambda(c,d)}", shift left=3, shorten <=4pt, shorten >=12pt, Rightarrow, from=0, to=3]
    \arrow["{\lambda(e,f)}"', shift right=3, shorten <=4pt, shorten >=12pt, Rightarrow, from=1, to=4]
    \arrow["{K(e,c)\otimes K(f,d)}"', shorten >=8pt, Rightarrow, from=3-3, to=2]
  \end{tikzcd}
\end{gather*}
Since vertical composition in $\SGame{}$ is given on losses by $+$, we can write this equation as
\begin{align}
  & \lambda(e\lenscirc c,f\lenscirc d) + K(e\otimes f,c\otimes d) \nonumber \\
  &\;= \lambda(e,f)\diamond \lambda(c,d) + K(e,c)\otimes K(f,d) \nonumber \\
  &\;= \lambda(e,f)_{c\otimes d} + \lambda(c,d)\circ(e'\otimes f')_{c\otimes d}  + K(e,c)\otimes K(f,d) \; . \label{eq:laxa-nat}
\end{align}
In each of the examples below, therefore, we establish the definition of the laxator $\lambda$ and check that it satisfies equation \ref{eq:laxa-nat}.

We will often use the notation $(-)_X$ to denote projection onto a factor $X$ of a monoidal product.

\subsubsection{Relative entropy}

\begin{prop} \label{prop:kl-lax-mon}
  The loss model $\KL$ of Proposition \ref{prop:kl-loss-model} is lax monoidal.
  Supposing that $(c,c'):(X,X)\lensto(Y,Y)$ and $(d,d'):(X',X')\lensto(Y',Y')$ are lenses in $\cat{B}$, the corresponding component $\lambda^{\KL}(c,d)$ of the laxator is given, for $\omega:I\klto X\otimes X'$ and $(y,y'):Y\otimes Y'$, by
  \[ \lambda^{\KL}(c,d)_\omega(y,y') := \E_{\substack{(x,x',m,m') \, \sim \\ (c'_{\omega_X}\otimes\, d'_{\omega_{X'}})(y,y')}} \left[ \log \frac{p_{\omega_X\otimes\omega_{X'}}(x,x')}{p_\omega(x,x')} \right] + \log \frac{p_{(c\otimes d)^{\smallground}\klcirc\omega}(y,y')}{p_{(c\otimes d)^{\smallground}\klcirc(\omega_X\otimes\omega_{X'})}(y,y')} \; . \]
\end{prop}

\subsubsection{Maximum likelihood estimation}

\begin{prop} \label{prop:mle-lax-mon}
  The loss model $\MLE$ of Proposition \ref{prop:mle-loss-model} is lax monoidal.
  Supposing that $(c,c'):(X,X)\lensto(Y,Y)$ and $(d,d'):(X',X')\lensto(Y',Y')$ are lenses in $\cat{B}$, the corresponding component $\lambda^{\MLE}(c,d)$ of the laxator is given, for $\omega:I\klto X\otimes X'$ and $(y,y'):Y\otimes Y'$, by
  \[ \lambda^{\MLE}(c,d)_\omega(y,y') := \log \frac{p_{(c\otimes d)^{\smallground}\klcirc(\omega_X\otimes\omega_{X'})}(y,y')}{p_{(c\otimes d)^{\smallground}\klcirc\omega}(y,y')} \; . \]
\end{prop}

\subsubsection{Free energy}

\begin{cor} \label{cor:fe-lax-mon}
  The loss model $\FE$ of Definition \ref{def:fe-loss-model} is lax monoidal.
  Supposing that $(c,c'):(X,X)\lensto(Y,Y)$ and $(d,d'):(X',X')\lensto(Y',Y')$ are lenses in $\cat{B}$, the corresponding component $\lambda^{\FE}(c,d)$ of the laxator is given, for $\omega:I\klto X\otimes X'$ and $(y,y'):Y\otimes Y'$, by
  \[ \lambda^{\FE}(c,d)_\omega(y,y') := \E_{(x,x')\sim(c'_{\omega_X}\otimes d'_{\omega_{X'}})(y,y')} \left[ \log \frac{p_{\omega_X\otimes\omega_{X'}}(x,x')}{p_\omega(x,x')} \right] \; . \]
\end{cor}

\subsubsection{Laplacian free energy}

\begin{prop} \label{prop:lfe-lax-mon}
  The loss model $\LFE$ of Propositions \ref{lemma:6-laplace-approx} and \ref{prop:lfe-loss-model} is lax monoidal.
  Supposing that $(c,c'):(X,X)\lensto(Y,Y)$ and $(d,d'):(X',X')\lensto(Y',Y')$ are lenses in $\cat{B}$, the corresponding component $\lambda^{\LFE}(c,d)$ of the laxator is given, for $\omega:I\klto X\otimes X'$ and $(y,y'):Y\otimes Y'$, by
  \[ \lambda^{\LFE}(c,d)_\omega(y,y') := \log \frac{p_{\omega_X\otimes\omega_{X'}}(\mu_{(c\otimes d)'_\omega}(y,y')_{XX'})}{p_\omega(\mu_{(c\otimes d)'_\omega}(y,y')_{XX'})} \]
  where $\mu_{(c\otimes d)'_\omega}(y,y')_{XX'}$ is the $(X\otimes X')$-mean of the Gaussian distribution $(c'_{\omega_X}\otimes d'_{\omega_{X'}})(y,y')$.
\end{prop}

\end{document}